\documentclass[]{interact}
\usepackage[caption=false]{subfig}

\usepackage[numbers,sort&compress]{natbib}
\bibpunct[, ]{[}{]}{,}{n}{,}{,}
\renewcommand\bibfont{\fontsize{10}{12}\selectfont}
 
\theoremstyle{plain}
\newtheorem{theorem}{Theorem}[section]
\newtheorem{lemma}[theorem]{Lemma}

\theoremstyle{definition}
\newtheorem{definition}[theorem]{Definition}

\theoremstyle{remark}
\newtheorem{remark}{Remark}

\usepackage{graphicx}
\usepackage[dvipsnames]{xcolor}

\usepackage{algorithm}
\usepackage{algorithmic}
\usepackage{caption}

\usepackage{makecell}
\usepackage{mathtools}
\usepackage{amsfonts}
\usepackage{amssymb}
\usepackage{url}
\usepackage{dirtytalk}
\usepackage{tikz}
\usepackage{nicefrac}
\usepackage{hyperref}
\usepackage{multirow}
\usepackage{cancel}

\begin{document}

\articletype{ARTICLE TEMPLATE}

\title{On Some Versions of Subspace Optimization Methods with Inexact Gradient Information \footnote{The research by F. Stonyakin in Section 4 has been financially supported by The Analytical Center for the Government of the Russian Federation (Agreement No. 70-2021-00143 01.11.2021, IGK 000000D730324P540002). The research by I.Kuruzov in Section 3 was supported by the Ministry of Economic Development of the Russian Federation (agreement No. 139-10-2025-034 dd. 19.06.2025, IGK 000000C313925P4D0002).
}
}

\author{
\name{Ilya Kuruzov\textsuperscript{a,b,c}\thanks{CONTACT Ilya Kuruzov. Email: kuruzov.ia@phystech.edu } and Fedor Stonyakin\textsuperscript{a,d}}
\affil{
\textsuperscript{a}Moscow Institute of Physics and Technology, Dolgoprudny, Russia;
\textsuperscript{b}Research Center for Artificial Intelligence,  Innopolis University, Innopolis, Russia;
\textsuperscript{c}National Research University Higher School of Economics, Moscow, Russia
\textsuperscript{d}V. Vernadsky Crimean Federal University, Simferopol, Republic of Crimea;
}
}

\maketitle

\begin{abstract}
It is well-known that accelerated gradient methods possess optimal complexity estimates for the class of convex smooth minimization problems. In many practical situations, it makes sense to work with inexact gradients. However, this can lead to the accumulation of corresponding inexactness in the theoretical estimates of the rate of convergence. We propose some modifications of first-order methods for convex optimization with an inexact gradient based on subspace optimization, such as Nemirovski's Conjugate Gradient method and the Sequential Subspace Optimization method. We study their convergence under different conditions on the inexactness both in the gradient value and in the accuracy of the solution of the subspace optimization subproblems. Besides this, we investigate a generalization of these results to the class of quasar-convex (weakly-quasi-convex) functions.
\end{abstract}


\section*{Introduction}\label{sec1_introduction}

First-order methods are an important class of approaches for optimization problems. They have different advantages: simple implementation, usually low cost of iterations and high performance for a wide class of functions (\cite{beck_first_ord}, \cite{gasnikov2017universal}, \cite{Polyak1987}). Nevertheless, there are many settings where a method has access only to an inexact gradient: gradient-free optimization in infinite-dimensional spaces \cite{Vasiliev}, inverse problems \cite{Kabanikhin}, saddle-point problems \cite{lan_first_ord}, \cite{inexact_primal_dual}. Therefore, methods with an inexact gradient are of interest to many researchers  \cite{DevolderThesis}, \cite{DevPaper}, \cite{Aspremont2008}, \cite{Polyak1987}, \cite{Vasin2021}.

Further, note that there are well-known results about the convergence of first-order methods and the optimality of accelerated methods for the class of convex functions \cite{d2021acceleration}, \cite{bubeck2015convex}, \cite{Nemirovski}, \cite{barre2020complexity}, \cite{NemirovskyYudin}. On the other hand, non-convex optimization appears in many practical problems. In particular, interest in such problems is growing because of deep learning \cite{lan_first_ord}, \cite{goodfellow}. One possible relaxation of convexity is quasar-convexity (also known as weak quasi-convexity). This class and non-convex examples are described in  \cite{AccWQC, Hardt, NearOpt}. {A recent work \cite{hermant} studies Nesterov Accelerated Gradient (NAG) for the class of strongly quasar-convex functions and shows, via a high-resolution ODE analysis, that acceleration can be achieved under a weaker curvature assumption. In contrast to that work, which relies on exact gradients, we focus on subspace optimization methods (SESOP and the conjugate gradient method) that operate with inexact gradient information.}

This paper continues the study of first-order methods based on subspace optimization. Such methods were considered in \cite{AccWQC, kuruzov_stonyakin}. One of them is the Sequential Subspace Optimization (SESOP) method \cite{SESOP_2005}. This method sequentially minimizes the objective over low-dimensional subspaces and converges to the solution. Recently, several interesting properties were demonstrated for this method. In particular, \cite{AccWQC} contains a proof of convergence of SESOP in the quasar-convex case. It demonstrates that this method converges with a rate similar to accelerated rates. Besides, \cite{kuruzov_stonyakin} proves its convergence in the case of an inexact gradient. Moreover, it shows that this method does not accumulate the error, in contrast to other known accelerated methods. However, there were no works devoted to the complexity of such methods in terms of the number of gradient calculations rather than the number of iterations. In this work, we propose to use the ellipsoid method \cite{ellipsoids, Kuruzov_Gladin} for the auxiliary low-dimensional problems. We obtain the total complexity (in terms of gradient calculations both for the outer iterations and for the subproblems) under the condition of convexity of the subproblems and in the presence of inexactness caused both by the inexact gradient and by the inexact solution of the subproblems. Let us stress that, although the convergence of SESOP with an inexact gradient and inexactly solved subproblems was analyzed in \cite{kuruzov_stonyakin}, the results of Section \ref{section:subproblems} are not a mere combination of known facts: we additionally establish an explicit and runtime-verifiable link between the accuracy of the subproblem solution in the function value and the inexactness conditions of \cite{kuruzov_stonyakin}, which makes it possible to estimate the total number of gradient calculations.

Besides SESOP, we consider a generalization of Nemirovski's Conjugate Gradient (CG) method. In Section \ref{section:conj_grad} we study its convergence for quasar-convex functions that satisfy the quadratic growth condition. Here and below, $\varepsilon$ denotes the desired accuracy in the function value, $\varepsilon_0$ is the initial residual in the function value, $\mu$ is the constant from the quadratic growth (or Polyak--{\L}ojasiewicz) condition, $L$ is the Lipschitz constant of the gradient, $\gamma$ is the quasar-convexity parameter and $\delta_1$ is the level of the additive inexactness in the gradient; see Section \ref{sec:problem_statement} for precise definitions. For the CG method we also demonstrate that it does not accumulate the additive gradient inexactness. Nevertheless, note that the quality of the solution obtained by the CG method can degrade to $O\left(\frac{1}{\gamma}\sqrt{\frac{\varepsilon_0}{\mu}}\delta_1\right)$. Also, we prove that to obtain the quality $\varepsilon$, the restarted version of the method needs $O\left(\sqrt{\frac{L}{\mu}}\log\frac{\varepsilon_0}{\varepsilon}\right)$ iterations for sufficiently small $\delta_1$. The main disadvantage of the result in this section is the requirement that the norm of the (inexact) gradient be sufficiently large at each iteration. To deal with it, we introduce a stopping condition that guarantees a compromise between the quality of the solution and the complexity of the algorithm.

Our contributions are the following:

\textbf{1. Linear convergence of the restarted CG method with an inexact gradient in the non-convex (quasar-convex) case.} We generalize the proof of convergence of Nemirovski's Conjugate Gradient method to the case of an inexact gradient. Moreover, we propose a stopping rule to reach the required quality.

\textbf{2. Complexity of auxiliary problems for SESOP and CG methods in the case of convex subproblems.} It is natural to solve the low-dimensional subproblems arising in these methods by methods for low-dimensional optimization. We consider the Ellipsoid Method and Multidimensional Dichotomy for these problems and estimate the total complexity when the subproblems are convex.

\section{Problem Statement}\label{sec:problem_statement}
Throughout the paper, we consider the unconstrained minimization problem
\begin{equation}\label{main_problem}
\min_{x\in\mathbb{R}^n} f(x),
\end{equation}
where $f$ is an $L$-smooth function ($\|\cdot\|$ is the usual Euclidean norm), i.e.
\begin{equation}\label{Lsmothness}
\|\nabla f(x) - \nabla f(y)\| \leqslant L \|x - y\| \quad \forall x, y \in \mathbb{R}^n,    
\end{equation}
and the methods have access only to an inexact gradient $g: \mathbb{R}^n \rightarrow \mathbb{R}^n$:
\begin{equation}\label{InexactGrad}
\|g(x) - \nabla f(x)\| \leqslant \delta_1 \quad \forall x\in\mathbb{R}^n,
\end{equation}
where $L>0$ and $\delta_1 \geq 0$. Below, $x^*$ denotes one of the exact solutions of problem \eqref{main_problem} and $f^* = f(x^*)$ denotes the optimal value.

The convergence results (in terms of the number of iterations) are formulated for quasar-convex problems.
\begin{definition}
Assume that $\gamma \in (0,1]$ and let $x^*$ be a minimizer of the differentiable function $f:\mathbb{R}^n \rightarrow \mathbb{R}$. The function $f$ is $\gamma$-quasar-convex with respect to $x^*$ if for all $x\in\mathbb{R}^n$,
\begin{equation}
    f(x^*)\geq f(x) + \frac{1}{\gamma}\langle\nabla f(x), x^*-x\rangle.
\end{equation}
\end{definition}

This class generalizes convex functions. It is also known as the class of weakly-quasi-convex functions. In the case $\gamma=1$ it coincides with the well-known star-convexity. \cite{NearOpt} demonstrates that the function $f(x)=(x^2+1/8)^{1/6}$ is quasar-convex but not convex or star-convex. Besides, the work \cite{wang_wibisono} demonstrates that some problems of training general linear models are quasar-convex.

Also, in this work, we will consider generalizations of strong convexity. The first considered condition is PL-condition.

\begin{definition}
The differentiable function $f$ satisfies the Polyak-Łojasiewicz condition (for brevity, we write PL-condition) for some constant $\mu>0$:
\begin{equation}\label{PL_cond}
f(x)-f^{*}\leqslant\frac{1}{2\mu}\|\nabla f(x)\|^{2}\quad \forall\,x\in\mathbb{R}^{n}.
\end{equation}
\end{definition}

This condition was first introduced in optimization in \cite{Polyak_1963}. Recently (see \cite{Karimi, Belkin}), it was proved that many practical problems satisfy this condition. In particular, it holds for over-parameterized non-linear systems.

Moreover, we obtain results on the performance of the well-known Conjugate Gradient method under a weaker condition than the PL-condition, namely, the quadratic growth condition (see \cite{Karimi}).

\begin{definition}
The differentiable function $f$ satisfies the quadratic growth condition (for brevity, we write QC-condition) for some constant $\mu>0$:
\begin{equation}\label{QC_cond}
f(x)-f^{*}\geqslant\frac{\mu}{2}\|x-x^*\|^{2}\quad \forall\,x\in\mathbb{R}^{n},
\end{equation}
where $x^*$ denotes a point of the solution set $X^* = \operatorname{Arg\,min} f$ nearest to $x$, i.e. $\|x - x^*\| = \operatorname{dist}(x, X^*)$.
\end{definition}

It is known (see \cite{Karimi}, Appendix A) that the PL-condition with constant $\mu$ implies the quadratic growth condition \eqref{QC_cond} with the same constant $\mu$ (with $\|x-x^*\| = \operatorname{dist}(x, X^*)$ the distance to the solution set). If the solution is unique, $x^*$ is simply this solution; in the results of Sections \ref{section:conj_grad} and \ref{section:subproblems}, $x^*$ is a fixed minimizer with respect to which quasar-convexity is assumed, and the quadratic growth condition is used at the points of the trajectory in the form $f(x)-f^*\geq\frac{\mu}{2}\operatorname{dist}^2(x, X^*)$, which is what the estimates actually require. Let us also note that the class of non-convex functions satisfying both quasar-convexity and the PL-condition (and, consequently, the quadratic growth condition) is non-empty. For example, the function
\begin{equation}\label{example_pl_quasar}
f(x) = \sum\limits_{i=1}^n \left(x_i^2 + 3\sin^2 x_i\right)
\end{equation}
is non-convex, $L$-smooth with $L = 8$ and satisfies the PL-condition with $\mu = 1/32$ (see \cite{Karimi}); at the same time, it is $\gamma$-quasar-convex with respect to $x^*=0$ with $\gamma = 0.49$. The latter fact is proved analytically in Appendix \ref{appendix_quasar_example}. The same profile $\varphi(t) = t^2 + 3\sin^2 t$ underlies the non-convex test problem \eqref{a_composed_problem} of Section \ref{section:numer_exp}.

\section{Subspace Optimization Methods}

Let us present the SESOP (Sequential Subspace Optimization) method from \cite{AccWQC, SESOP_2005, kuruzov_stonyakin} (see Algorithm \ref{alg:sesop}). The first step of this method is the construction of a subspace for further optimization. At each iteration, there are three important directions: the gradient at the current point, the direction from the starting point to the current one, and a weighted sum of the gradients from all previous iterations. Note that all these directions can be calculated using only $x_k, x_0, g(x_k)$ and the previously accumulated weighted sum. In other words, this method does not require too much additional memory in comparison with other first-order methods like Gradient Descent.

\begin{algorithm}[ht]
	\caption{A modification of the SESOP method with an inexact gradient}
  \label{alg:sesop}
  \begin{algorithmic}[1]
  \REQUIRE objective function $f$ with an inexact gradient $g$, initial point $x_0$, number of iterations $T$.
\STATE $w_0 = 1$
  \FOR{$k =0, \hdots, T-1$}
      \STATE Construct Subspace:
      $\mathbf{d}_k^0 = g(x_k), \quad \mathbf{d}_k^1 = x_k-x_0,\quad\mathbf{d}_k^2 = \sum\limits_{i=0}^k w_i g(x_i).$

      \STATE Find the optimal step
      \begin{equation}\label{eqsubproblem_sesop}
          \tau_k \leftarrow \arg\min_{\tau\in\mathbb{R}^3}f\left(x_k + \sum\limits_{i=1}^3 \tau_i \mathbf{d}_k^{i-1} \right)
      \end{equation}
      \STATE $x_{k+1} \leftarrow x_k + \sum\limits_{i=1}^3 (\tau_k)_i \mathbf{d}_k^{i-1}$
      \STATE Update $w:$ $w_{k+1} = \frac{1}{2}+\sqrt{\frac{1}{4}+w_{k}^2}$
  \ENDFOR
\RETURN $x_T$
	\end{algorithmic}
\end{algorithm}

The next step is the optimization over the three-dimensional subspace. It is the most complex step. Note that even if the original problem is quasar-convex, the auxiliary problem may not have good properties. It is the bottleneck of this method, and it will be discussed in Section \ref{section:subproblems} for the case of convex subproblems.

The final steps are the calculation of the new point and the update of the weight for the direction of the weighted gradient sum. These steps do not require additional computations. The convergence rate (per iteration) of this method is presented in \cite{SESOP_2005}. This result was recently generalized to the quasar-convex case in \cite{AccWQC}. Further, it was proved in \cite{kuruzov_stonyakin} that this method is robust to inexactness in the gradient.

\begin{algorithm}[ht]
	\caption{A modification of Nemirovski's Conjugate Gradient method with an inexact gradient}
  \label{alg:CG}
  \begin{algorithmic}[1]
  \REQUIRE objective function $f$ with an inexact gradient $g$, initial point $x_0$, number of iterations $T$.

  \STATE $\mathbf{q}_0 = 0$
  
  \FOR{$k =0, \hdots, T-1$}
      \STATE Solve 2-dimensional problem
      \begin{equation}\label{eqsubproblem_cg}
          \hat{x}_k \leftarrow \arg\min_{x \in \mathcal{X}_k}f(x), \quad\text{where } \mathcal{X}_k=x_0 +\text{Lin}(x_k-x_0, \mathbf{q}_k)
      \end{equation}

      \STATE Make gradient step: $x_{k+1} = \hat{x}_k - \frac{1}{2L}g(\hat{x}_k)$
      \STATE $\mathbf{q}_{k+1} = \mathbf{q}_{k} + g(\hat{x}_k)$
  \ENDFOR
\RETURN $x_T$
	\end{algorithmic}
\end{algorithm}

Another method considered in \cite{AccWQC} is Nemirovski's Conjugate Gradient method (see Algorithm \ref{alg:CG}). It is a well-known method with rather high performance in practice. Besides, it is known that conjugate gradient methods have a closed form for quadratic minimization problems and are optimal for them. There are different ways to generalize such methods to non-quadratic problems. In this paper, we consider Nemirovski's version \cite{Nemirovski}. In \cite{Nemirovski}, the theoretical convergence rate of the CG method was a consequence of the following properties: 1) smoothness of the function, 2) strong convexity of the function, 3) orthogonality of the gradient at the current point to the direction from the starting point to the current one, 4) orthogonality of the gradient at the current point to the sum of the gradients from all previous iterations.

In \cite{AccWQC}, it was proved that this method converges if the strong convexity condition is replaced by quasar-convexity and the quadratic growth condition. The last two properties (the orthogonality conditions) are consequences of the optimization over the 2-dimensional subspace (see step 2 in Algorithm \ref{alg:CG}). Nevertheless, this orthogonality becomes inexact in the case of an inexact solution of the auxiliary subproblem. Sections \ref{section:conj_grad} and \ref{section:subproblems} are devoted to this problem.

\section{Convergence of CG Method}\label{section:conj_grad}

In \cite{AccWQC}, the authors obtained a result on the convergence rate of Nemirovski's Conjugate Gradient method. In this work, we show that Algorithm \ref{alg:CG} can also work with an additively inexact gradient when the inexactness is not too large.

In this section, we consider two versions of this approach. Firstly, we consider the classic Nemirovski's Conjugate Gradient method. However, without restarts this method does not attain the linear convergence rate in the strongly convex (or quadratic growth) case. The restart technique allows one to overcome this issue. So, after the analysis of the classic approach, we will consider its work with the restart technique. Let us stress that Theorem \ref{CG_theorem} below is not a final convergence result by itself: it only quantifies the guaranteed progress of Algorithm \ref{alg:CG} after $T$ iterations and serves as an intermediate step towards the restarted scheme, for which the convergence rate is given in Theorem \ref{CG_theorem_restarts}.
\subsection{Convergence with Inexact Gradient}
\label{CG_section1}

We start from the classic Nemirovski's CG method. To provide the analysis under an inexact gradient, we need the following auxiliary lemma. In the notation of Algorithm \ref{alg:CG}, $q_T = \sum_{k=0}^{T-1} g(\hat{x}_k)$.

\begin{lemma}
\label{CG_lemma}
Let the objective function $f$ be $L$-smooth and $\gamma$-quasar-convex with respect to $x^*$. Also for the inexact gradient $g:\mathbb{R}^n\rightarrow\mathbb{R}^n$ there is some constant $\delta_1\geq 0$ such that for all $x\in\mathbb{R}^n$:     \begin{equation}
    \label{g_cond_i}
        \|g(x)-\nabla f(x)\|\leq\delta_1.
    \end{equation}
    Let, moreover, the auxiliary two-dimensional subproblems \eqref{eqsubproblem_cg} in Algorithm \ref{alg:CG} be solved exactly. Then the following inequality holds:
\begin{equation}
\label{q_t_est}
\|q_T\|\leq 4\delta_1 T +\left(\sum\limits_{k=0}^{T-1} \|g(\hat{x}_k)\|^2\right)^{\frac{1}{2}}.
\end{equation}
\end{lemma}
The proof of Lemma \ref{CG_lemma} is given in Appendix \ref{appendix_lemma_proof}.

Using Lemma \ref{CG_lemma}, we can generalize the result of Theorem 2 from \cite{AccWQC} to the case of an inexact gradient. Finally, we have the following result.

\begin{theorem}
\label{CG_theorem}
Let the objective function $f$ be $L$-smooth and $\gamma$-quasar-convex with respect to $x^*$. Also for the inexact gradient $g:\mathbb{R}^n\rightarrow\mathbb{R}^n$ there is some constant $\delta_1\geq 0$ such that for all $x\in\mathbb{R}^n$ $\|g(x)-\nabla f(x)\|\leq\delta_1.$
Moreover, let the function satisfy the quadratic growth condition $f(x)-f^*\geq \frac{\mu}{2}\|x-x^*\|^2$ and let the auxiliary subproblems \eqref{eqsubproblem_cg} be solved exactly. Then, if at all iterations $k = 0, \ldots, T-1$ the condition $\|g(\hat{x}_k)\|\geq 2\delta_1$ holds, Algorithm \ref{alg:CG} obtains a point $x_T$ such that
\begin{equation}
    \label{est_CG}
     f(x_T)-f^*\leq \beta\epsilon_0 +\frac{7}{\gamma}\sqrt{\frac{2\varepsilon_0}{\mu}}\delta_1
\end{equation}
after
$$T=\left\lceil\frac{2}{\gamma \beta} \sqrt{\frac{2(1-\beta)L}{\mu}}\right\rceil+1$$
iterations, where $\epsilon_0=f(x_0) - f^*$ and $\beta\in(0,1)$ is an arbitrary constant parameter.
\end{theorem}

\begin{proof}
Assume, on the contrary, that
\begin{equation}
\label{contr_assump}
\varepsilon_T := f(x_T)-f^*\geq \beta \varepsilon_0 + c\delta_1,
\end{equation}
where $\varepsilon_k := f(x_k) - f^*$ and $c>0$ is a constant that will be chosen later. We will show that for $T$ from the statement of the theorem this assumption leads to a contradiction.

{\it Step 1 (descent inequality and its consequences).} Since $x_{k+1} = \hat{x}_k - \frac{1}{2L}g(\hat{x}_k)$, the $L$-smoothness of $f$ implies
$$f(x_{k+1}) \leq f(\hat{x}_k) - \frac{1}{2L}\langle \nabla f(\hat{x}_k), g(\hat{x}_k)\rangle + \frac{1}{8L}\|g(\hat{x}_k)\|^2.$$
Using \eqref{g_cond_i} and the Cauchy--Schwarz inequality, $\langle \nabla f(\hat{x}_k), g(\hat{x}_k)\rangle \geq \|g(\hat{x}_k)\|^2 - \delta_1\|g(\hat{x}_k)\|$, whence
\begin{equation}
\label{descent_ineq}
f(x_{k+1}) \leq f(\hat{x}_k) - \frac{3}{8L}\|g(\hat{x}_k)\|^2 + \frac{\delta_1}{2L}\|g(\hat{x}_k)\| \leq f(\hat{x}_k) - \frac{1}{4L}\|g(\hat{x}_k)\|^2 + \frac{\delta_1^2}{2L},
\end{equation}
where in the last step we applied Young's inequality $\frac{\delta_1}{2L}\|g(\hat{x}_k)\| \leq \frac{1}{8L}\|g(\hat{x}_k)\|^2 + \frac{\delta_1^2}{2L}$. Therefore,
\begin{equation}
\label{L_ineq}
\|g(\hat{x}_k)\|^2 \leq 4L\left(f(\hat{x}_k)-f(x_{k+1})\right)+2\delta_1^2.
\end{equation}
Moreover, since $\hat{x}_k$ is an exact minimizer of $f$ over the affine subspace $\mathcal{X}_k \ni x_k$, we have $f(\hat{x}_k)\leq f(x_k)$, so \eqref{L_ineq} yields
\begin{equation}
\label{L_ineq2}
\|g(\hat{x}_k)\|^2 \leq 4L\left(\varepsilon_k-\varepsilon_{k+1}\right)+2\delta_1^2.
\end{equation}
Summing \eqref{L_ineq2} over $k = 0, \ldots, T-1$ and using $\varepsilon_0 - \varepsilon_T \leq (1-\beta)\varepsilon_0$ (which follows from \eqref{contr_assump}), we obtain
\begin{equation}
\label{f_k_est_eps}
\sum\limits_{k=0}^{T-1}\|g(\hat{x}_k)\|^2 \leq 4L(\varepsilon_0 - \varepsilon_T) + 2T\delta_1^2 \leq 4L(1-\beta)\varepsilon_0 + 2T\delta_1^2.
\end{equation}
Note also that the assumption $\|g(\hat{x}_k)\| \geq 2\delta_1$ together with the middle expression of \eqref{descent_ineq} gives $f(x_{k+1}) \leq f(\hat{x}_k) - \frac{\delta_1^2}{2L} \leq f(\hat{x}_k)$. Combining this with $f(\hat{x}_k)\leq f(x_k)$, we conclude that the sequence $f(x_0) \geq f(\hat{x}_0)\geq f(x_1) \geq f(\hat{x}_1) \geq \ldots \geq f(x_T)$ is non-increasing. In particular, for every $k \leq T-1$,
\begin{equation}
\label{monotone_chain}
f(\hat{x}_k) - f^* \geq \varepsilon_T \geq \beta\varepsilon_0 + c\delta_1,
\end{equation}
i.e. assumption \eqref{contr_assump}, made for the final point, automatically extends to all intermediate points $\hat{x}_k$.

{\it Step 2 (quasar-convexity).} From quasar-convexity we have the estimate
$$f(\hat{x}_k)-f^* \leq \frac{1}{\gamma} \langle \nabla f(\hat{x}_k), \hat{x}_k - x^*\rangle.$$
Since $\hat{x}_k$ is an exact minimizer of $f$ over the affine subspace $\mathcal{X}_k = x_0 + \mathrm{Lin}(x_k - x_0, q_k)$, the gradient $\nabla f(\hat{x}_k)$ is orthogonal to the directions of this subspace; in particular, as $\hat{x}_k - x_0 \in \mathrm{Lin}(x_k - x_0, q_k)$, we have $\langle \nabla f(\hat{x}_k), \hat{x}_k - x_0\rangle = 0$. Hence
$$f(\hat{x}_k)-f^* \leq \frac{1}{\gamma} \langle \nabla f(\hat{x}_k), x_0 - x^*\rangle \leq \frac{1}{\gamma} \langle g(\hat{x}_k), x_0 - x^*\rangle + \delta_1 \frac{R}{\gamma},$$
where $R := \|x_0 - x^*\|$ and in the last step we used \eqref{g_cond_i} together with the Cauchy--Schwarz inequality. Combining this with \eqref{monotone_chain} and summing over $k = 0, \ldots, T-1$ (recall that $q_T = \sum_{k=0}^{T-1} g(\hat{x}_k)$), we get
$$T\beta \varepsilon_0 + cT\delta_1 \leq \frac{1}{\gamma} \langle q_T, x_0 - x^*\rangle + \delta_1 \frac{RT}{\gamma},$$
which, by the Cauchy--Schwarz inequality, implies
\begin{equation}
\label{est1}
\|q_T\|\|x_0-x^*\| \geq T \gamma \beta \varepsilon_0 + \delta_1T \left(c\gamma - R\right).
\end{equation}

{\it Step 3 (combining the bounds).} From Lemma \ref{CG_lemma} and inequality \eqref{f_k_est_eps}, using $\sqrt{a+b}\leq\sqrt{a}+\sqrt{b}$ and $\sqrt{2T}\leq \sqrt{2}\,T$ for $T\geq 1$, we have
\begin{equation}
\label{est2}
\|q_T\|\leq 4\delta_1 T + \sqrt{4L (1-\beta)\varepsilon_0 + 2T\delta_1^2} \leq 6\delta_1 T + 2\sqrt{L (1-\beta)\varepsilon_0}.
\end{equation}
On the other hand, applying the quadratic growth condition at the point $x_0$, we obtain the following estimate for $\|x_0-x^*\|$:
\begin{equation}
\label{est3}
\|x_0-x^*\|\leq \sqrt{\frac{2\varepsilon_0}{\mu}}.
\end{equation}
Here $x^*$ is taken to be a minimizer nearest to $x_0$, with respect to which the quasar-convexity is assumed (in the case of a unique solution this is automatic); this is exactly the setting in which the quadratic growth condition \eqref{QC_cond} yields \eqref{est3}.
Combining inequalities \eqref{est1}--\eqref{est3}, we obtain
\begin{equation}
\label{T_est}
\left(6\delta_1 T + 2\sqrt{L (1-\beta)\varepsilon_0}\right)\sqrt{\frac{2\varepsilon_0}{\mu}}\geq T \gamma \beta \varepsilon_0+\delta_1T \left(c\gamma - R \right),
\end{equation}
which can be rewritten in the following form:
$$T \gamma \beta \varepsilon_0 + \delta_1 T\left(c\gamma - R -6\sqrt{\frac{2\varepsilon_0}{\mu}} \right)\leq 2\varepsilon_0 \sqrt{\frac{2(1-\beta)L}{\mu}}.$$
Choosing $c = \frac{1}{\gamma}\left(R +6\sqrt{\frac{2\varepsilon_0}{\mu}}\right)$, the second term on the left-hand side vanishes and we obtain the following estimate on $T$:
$$T\leq \frac{2}{\gamma \beta} \sqrt{\frac{2(1-\beta)L}{\mu}}.$$
Thus, for $T=\left\lceil\frac{2}{\gamma \beta} \sqrt{\frac{2(1-\beta)L}{\mu}}\right\rceil + 1 > \frac{2}{\gamma \beta} \sqrt{\frac{2(1-\beta)L}{\mu}}$, assumption \eqref{contr_assump} leads to a contradiction. Consequently, for such $T$ we have $\varepsilon_T\leq \beta \varepsilon_0 + c\delta_1.$ Finally, by \eqref{est3} we have $R \leq \sqrt{\frac{2\varepsilon_0}{\mu}}$, so $c \leq \frac{7}{\gamma}\sqrt{\frac{2\varepsilon_0}{\mu}}$ and after $T$ iterations we obtain a point $x_T$ that satisfies the following estimate:
$$f(x_T)-f^* \leq \beta \varepsilon_0 + \frac{7}{\gamma}\sqrt{\frac{2\varepsilon_0}{\mu}}\delta_1.$$
\end{proof}

\begin{remark}
Note that the quadratic growth condition is met when the objective function $f$ satisfies the well-known PL-condition \eqref{PL_cond}: this implication (with the same constant $\mu$) is proved in \cite{Karimi} (Appendix A). It is well known \cite{nesterov2006cubic, Karimi, gasnikov2017universal} that under additional smoothness assumptions, standard non-accelerated iterative methods for such functions (Gradient Descent, Cubic Regularized Newton method, etc.) converge as if $f$ were a $\mu$-strongly convex function. For accelerated methods, such results are not known. So we were motivated to find additional sufficient conditions that guarantee convergence for properly chosen accelerated methods. In this section we observe that such a condition could be the $\gamma$-quasar-convexity of $f$.
\end{remark}

\begin{remark}
Note that if the function $f$ meets the PL-condition \eqref{PL_cond} and we stop our method when $\|g(x_k)\|\leq 2\delta_1$, then $\|\nabla f(x_k)\| \leq \|g(x_k)\| + \delta_1 \leq 3\delta_1$ and, by \eqref{PL_cond}, $f(x_k)- f^*\leq \frac{9\delta_1^2}{2\mu}$. Note that in \cite{stonyakin} the authors proved that there are no first-order methods with a $\delta_1$-inexact gradient that can converge better than $O\left(\frac{\delta_1^2}{\mu}\right)$ in the general case.
\end{remark}

\subsection{Nemirovski's CG Method with Restarts}
\label{CG_section_restarts}

It is well known that the restart technique can significantly improve the convergence of conjugate gradient methods. To use it, we run the algorithm for $T$ iterations, and after that start the same algorithm from the obtained point (see Algorithm \ref{alg:CG_restarts}).
\begin{algorithm}
	\caption{Restarted Nemirovski's Conjugate Gradient method with an inexact gradient}
  \label{alg:CG_restarts}
  \begin{algorithmic}[1]
  \REQUIRE objective function $f$ with an inexact gradient $g$, initial point $x_0$, number of iterations $T$, number of restarts $K$.

  \STATE $\mathbf{q}_0 = 0$
  
  \FOR{$k =1, \hdots, K$}
      \STATE Run Algorithm \ref{alg:CG} for start point $x_{k-1}$ and $T$ iterations:
      $$x_k = \text{CG}(f,g,x_{k-1}, T)$$
  \ENDFOR
\RETURN $x_K$
	\end{algorithmic}
\end{algorithm}

Similarly to \cite{AccWQC}, we can obtain the following result.

\begin{theorem}
\label{CG_theorem_restarts}
Let the objective function $f$ be $L$-smooth and $\gamma$-quasar-convex with respect to $x^*$. Also for the inexact gradient $g:\mathbb{R}^n\rightarrow\mathbb{R}^n$ there is some constant $\delta_1\geq 0$ such that for all $x\in\mathbb{R}^n$ $\|g(x)-\nabla f(x)\|\leq\delta_1$ and $\delta_1^2 \leq  \frac{\gamma^2 \alpha^2 \mu \varepsilon}{98}$ for some $\alpha\in (0,1)$, where $\varepsilon$ is the desired accuracy.
Moreover, let the function satisfy the quadratic growth condition $f(x)-f^*\geq \frac{\mu}{2}\|x-x^*\|^2$ and let the auxiliary subproblems \eqref{eqsubproblem_cg} be solved exactly. Then, if at all inner iterations the condition $\|g(\hat{x}_k)\|\geq 2\delta_1$ is met, the restarted CG method (Algorithm \ref{alg:CG_restarts}) obtains an output point $\hat{x} = x_K$ such that
\begin{equation}
     f(\hat{x})-f^*\leq \varepsilon
\end{equation}
after
$$K = \left\lceil\frac{2}{1-\alpha}\log\frac{\varepsilon_0}{\varepsilon}\right\rceil$$
restarts with
$$T=\left\lceil\frac{4}{\gamma} \sqrt{\frac{L}{\mu}}\frac{\sqrt{1+\alpha}}{1-\alpha}\right\rceil+1$$
iterations per restart, where $\epsilon_0=f(x_0) - f^*$.
\end{theorem}
\begin{proof}
Fix $\alpha \in (0,1)$ and set $\beta = \frac{1-\alpha}{2}$ in Theorem \ref{CG_theorem}. Consider one run of Algorithm \ref{alg:CG} (one restart) started from a point with the residual $\varepsilon_0'$ in the function value. By Theorem \ref{CG_theorem}, after
$$T=\left\lceil\frac{2}{\gamma \beta} \sqrt{\frac{2(1-\beta)L}{\mu}}\right\rceil + 1 = \left\lceil\frac{4}{\gamma(1-\alpha)} \sqrt{\frac{(1+\alpha)L}{\mu}}\right\rceil + 1$$
iterations, we obtain a point with the residual at most $\beta \varepsilon_0' + \frac{7}{\gamma}\sqrt{\frac{2\varepsilon_0'}{\mu}}\delta_1.$ If
\begin{equation}
\label{CG_delta_condition}
\frac{7}{\gamma}\sqrt{\frac{2\varepsilon_0'}{\mu}}\delta_1\leq \alpha\varepsilon_0',
\end{equation}
then the residual after the run is at most $\tilde{\beta}\varepsilon_0'$, where $\tilde{\beta}=\beta + \alpha = \frac{1+\alpha}{2} < 1$. Note that condition \eqref{CG_delta_condition} can be rewritten in the following form:
\begin{equation}
\label{CG_delta_condition1}
\delta_1^2 \leq  \frac{\gamma^2 \alpha^2 \mu \varepsilon_0'}{98}.
\end{equation}
Let $\varepsilon_j$ denote the residual after the $j$-th restart. As shown in the proof of Theorem \ref{CG_theorem}, the function values do not increase during the work of the method. Hence, as long as $\varepsilon_j \geq \varepsilon$, the assumption
\begin{equation}
\label{CG_delta_condition2}
\delta_1^2 \leq  \frac{\gamma^2 \alpha^2 \mu \varepsilon}{98}
\end{equation}
of the theorem implies \eqref{CG_delta_condition1} for $\varepsilon_0' = \varepsilon_j$, and therefore each restart contracts the residual: $\varepsilon_{j+1} \leq \tilde{\beta}\varepsilon_j$. Consequently, to guarantee $\varepsilon_K \leq \tilde{\beta}^K \varepsilon_0 \leq \varepsilon$, it suffices to make
$$K = \left\lceil\frac{\log({\varepsilon_0}/{\varepsilon})}{\log({1}/{\tilde{\beta}})}\right\rceil \leq  \left\lceil\frac{2}{1-\alpha}\log\frac{\varepsilon_0}{\varepsilon}\right\rceil$$
restarts, where we used $\log\frac{1}{\tilde{\beta}} = -\log\frac{1+\alpha}{2} \geq 1 - \frac{1+\alpha}{2} = \frac{1-\alpha}{2}$. Thus, after $K$ restarts the method obtains a point $x_{K}$ such that $f(x_{K}) - f^* \leq \varepsilon.$
\end{proof}

\begin{remark}
Generally, without condition \eqref{CG_delta_condition2}, after $K$ restarts (i.e. after $KT$ iterations in total) we obtain a point $x_K$ such that
\begin{equation*}
     f(x_K)-f^*\leq \beta^K \varepsilon_0 + \left(\sum\limits_{j=0}^{K-1}\beta^j\right)\frac{7}{\gamma}\sqrt{\frac{2\epsilon_0}{\mu}}\delta_1
     \leq \beta^K \epsilon_0 + \frac{7}{\gamma (1-\beta)}\sqrt{\frac{2\epsilon_0}{\mu}}\delta_1.
\end{equation*}
Here we can see that it cannot be guaranteed that Nemirovski's Conjugate Gradient method converges to a quality better than $O\left(\frac{1}{\gamma}\sqrt{\frac{\varepsilon_0}{\mu}}\delta_1\right).$
\end{remark}

\begin{remark}
Under the assumptions of Theorem \ref{CG_theorem_restarts}, Algorithm \ref{alg:CG_restarts} requires the total number of gradient computations $KT = O\left(\frac{1}{\gamma}\sqrt{\frac{L}{\mu}}\log\frac{\varepsilon_0}{\varepsilon}\right)$ (for a fixed $\alpha$) to approach the quality $\varepsilon$.
\end{remark}

As we mentioned above, first-order methods cannot approach a quality better than $O(\frac{\delta_1^2}{\mu})$ for strongly convex functions. Consequently, it is true for functions that meet the PL-condition \eqref{PL_cond} or the quadratic growth condition \eqref{QC_cond}. So, let us consider an accuracy of the order $O\left(\frac{\delta_1^2}{\gamma^2\mu}\right)$ acceptable for the function value and {\it agree to terminate Algorithm \ref{alg:CG_restarts} if the condition $\|g(\hat{x}_k)\|\leq \frac{8}{\gamma}\delta_1$ is satisfied; the criterion is checked at every point $\hat{x}_k$ where the inexact gradient is computed.}

In such a case, we have the following alternatives. On the one hand, if all the computed gradients are large enough, then the algorithm obtains a point with a sufficiently good quality after a specific number of restarts and iterations. On the other hand, if there is a point with a small enough norm of the inexact gradient, then this point already has a good quality with respect to the function value because of the PL-condition. Finally, this gives us the following theoretical result about the work of the CG method with the stopping rule.

\begin{theorem}\label{thm-main}
Let the objective function $f$ be $L$-smooth and $\gamma$-quasar-convex with respect to $x^*$. Also for the inexact gradient $g:\mathbb{R}^n\rightarrow\mathbb{R}^n$ there is some constant $\delta_1\geq 0$ such that for all $x\in\mathbb{R}^n$ $\|g(x)-\nabla f(x)\|\leq\delta_1$.
Moreover, let the function satisfy the PL-condition \eqref{PL_cond}, let the auxiliary subproblems be solved exactly, and let $\varepsilon=\frac{196}{\gamma^2 \mu} \delta_1^2$ and $\alpha\in\left[\frac{1}{\sqrt{2}},1\right)$.

Let one of the following alternatives hold:
\begin{enumerate}
    \item The stopping criterion is never satisfied and the restarted Nemirovski's Conjugate Gradient method (Algorithm \ref{alg:CG_restarts}) makes $$K = \left\lceil\frac{2}{1-\alpha}\log\frac{\varepsilon_0}{\varepsilon}\right\rceil$$
restarts with
$$T=\left\lceil\frac{4}{\gamma} \sqrt{\frac{L}{\mu}}\frac{\sqrt{1+\alpha}}{1-\alpha}\right\rceil+1$$
iterations per each restart;
\item At some inner iteration, the stopping criterion $\|g(\hat{x}_N)\|\leq \frac{8}{\gamma}\delta_1$ is satisfied for the first time at a computed point $\hat{x}_N$.
\end{enumerate}
Then for the output point $\widehat{x}$ ($\widehat{x} = x_{K}$ in the first case and $\widehat{x} = \hat{x}_{N}$ in the second case), the following inequality holds:
$$
f(\widehat{x})-f^{*}\leqslant\frac{196\delta_1^{2}}{\gamma^2 \mu}.
$$
\end{theorem}

\begin{proof}
First, recall that the PL-condition \eqref{PL_cond} with constant $\mu$ implies the quadratic growth condition \eqref{QC_cond} with the same constant (see \cite{Karimi}, Appendix A). Note that the two alternatives are exhaustive: either the stopping criterion $\|g(\hat{x}_k)\|\leq\frac{8}{\gamma}\delta_1$ is satisfied at some computed point (the second alternative), or it is never satisfied, in which case $\|g(\hat{x}_k)\| > \frac{8}{\gamma}\delta_1 \geq 8\delta_1 > 2\delta_1$ at all computed points (recall $\gamma\leq 1$), which is exactly the precondition of Theorem \ref{CG_theorem_restarts} (the first alternative). Consider the second alternative. If $\|g(\hat{x}_N)\|\leq \frac{8}{\gamma}\delta_1$, then $\|\nabla f(\hat{x}_N)\| \leq \left(\frac{8}{\gamma}+1\right)\delta_1 \leq \frac{9}{\gamma}\delta_1$ (recall that $\gamma\leq 1$), and the PL-condition yields
$$f(\hat{x}_N) - f^* \leq \frac{1}{2\mu}\|\nabla f(\hat{x}_N)\|^2 \leq \frac{81\delta_1^2}{2\gamma^2\mu} \leq \frac{196\delta_1^2}{\gamma^2\mu}.$$
Now consider the first alternative. Since the stopping criterion was never satisfied, at all computed points we have $\|g(\hat{x}_k)\| > \frac{8}{\gamma}\delta_1 \geq 2\delta_1$. Moreover, for $\varepsilon=\frac{196}{\gamma^2 \mu} \delta_1^2$ and $\alpha \geq \frac{1}{\sqrt 2}$ we have
$$\frac{\gamma^2\alpha^2\mu\varepsilon}{98} = 2\alpha^2\delta_1^2 \geq \delta_1^2,$$
so condition \eqref{CG_delta_condition2} of Theorem \ref{CG_theorem_restarts} holds. Applying Theorem \ref{CG_theorem_restarts}, we obtain $f(x_K) - f^* \leq \varepsilon = \frac{196\delta_1^2}{\gamma^2\mu}$.
\end{proof}

We can see that the achieved accuracy $O\left(\frac{\delta_1^2}{\gamma^2\mu}\right)$ matches (up to a constant factor and the factor $\frac{1}{\gamma^2}$) the lower bound $\Omega\left(\frac{\delta_1^2}{\mu}\right)$ from \cite{stonyakin} on the accuracy attainable by first-order methods with a $\delta_1$-inexact gradient. As for the iteration complexity $O\left(\frac{1}{\gamma}\sqrt{\frac{L}{\mu}}\log\frac{\varepsilon_0}{\varepsilon}\right)$, to the best of our knowledge, lower bounds specific to the considered class of quasar-convex functions satisfying the PL (or quadratic growth) condition are not known. Note, however, that for $\gamma = 1$ the class contains all smooth strongly convex functions (which are $1$-quasar-convex and satisfy quadratic growth), so the classical lower bound $\Omega\left(\sqrt{\frac{L}{\mu}}\log\frac{1}{\varepsilon}\right)$ for smooth strongly convex minimization (see \cite{Nemirovski}, \cite{nesterov2018lectures}) applies to this subclass; our estimate coincides with this optimal rate for $\gamma = 1$. A matching lower bound for general $\gamma < 1$ (in particular, whether the factor $1/\gamma$ is necessary) remains an open question. Nevertheless, in this case we have an additional parameter for tuning --- the frequency of restarts.

\section{Auxiliary Low-Dimensional Subproblems}\label{section:subproblems}

In this section, we assume that all auxiliary subproblems in Algorithms \ref{alg:sesop} and \ref{alg:CG} are convex. Namely, subproblems in step 4 of Algorithm \ref{alg:sesop} and in step 2 of Algorithm \ref{alg:CG} are convex. Let us note that this is a restrictive assumption: for a general quasar-convex $f$, the restriction of $f$ to a low-dimensional affine subspace need not be convex. This assumption holds, for example, when $f$ itself is convex. We adopt it because the complexity results for the low-dimensional cutting-plane-type methods used below (the ellipsoid method and the multidimensional dichotomy) are available only in the convex case (see \cite{Kuruzov_Gladin}); in practice, the subproblems can be solved by any suitable local method (see also Section \ref{section:numer_exp}).

\subsection{Ellipsoid Methods for SESOP}
\label{sesop_int_problem}

To estimate the number of gradient calculations, we need to choose some procedure for the optimization over the subspace in step 4 of Algorithm \ref{alg:sesop}. It is a three-dimensional problem, so we can use methods for low-dimensional problems. Examples of such methods are the ellipsoid method (see \cite{ellipsoids}), Vaidya's method (see \cite{vaidya1989new}) and the dichotomy method for a hypercube (see \cite{Kuruzov_Gladin}). For all these methods, there are results on their work with an inexact gradient (see \cite{Kuruzov_Gladin}). The dichotomy method has a worse estimate for the number of calculations than the other methods. Nevertheless, it demonstrates rather good performance in the two-dimensional case. Therefore, we consider it for the CG method below. The ellipsoid and Vaidya's methods require $O\left(\log\frac{1}{\varepsilon}\right)$ gradient calculations. So, in the current work we chose the ellipsoid method (see Algorithm \ref{alg:ellispoids}) for the subproblem.

For the ellipsoid method, there is the following estimate (see Theorem 2 in \cite{Kuruzov_Gladin}). If Algorithm \ref{alg:ellispoids} is run on a ball $\mathcal{B}\subset \mathbb{R}^n$ of radius $\mathcal{R}$ in the $n$-dimensional space and the constant $B$ is such that $\max_{x\in\mathcal{B}} f(x) - \min_{x\in\mathcal{B}} f(x)\leq B$, then the ellipsoid method with a $\delta$-subgradient converges to the solution with the following speed:

\begin{equation}
    \label{ellispoid_estimate}
    f(x_N) - f(x^*)\leq B\exp\left(-\frac{N}{2n^2}\right)+\delta
\end{equation}

In our case, $n$ is equal to 3, the dimension of the subproblem. So, according to \eqref{ellispoid_estimate}, when $\delta\leq \frac{\varepsilon}{2}$, to approach the quality $\varepsilon$ we need
\begin{equation}
\label{N_ellispoids}
    N\geq 18 \ln \frac{2B}{\varepsilon}
\end{equation}
iterations of Algorithm \ref{alg:ellispoids}.

In this section, we will estimate the work of the SESOP algorithm in two settings:
\begin{itemize}
    \item one has an exact low-dimensional gradient (the gradient of the subproblem) but only an inexact gradient of the full problem;
    \item one has only inexact gradients both for the low-dimensional subproblem and for the full problem.
\end{itemize}

Obviously, if we have an inexact gradient of the full problem, then we are able to calculate an inexact gradient of the auxiliary problem too. On the other hand, the low-dimensional oracle can be significantly more exact in some cases because of a simpler calculation. So, let us start from considering the first setting.

In the paper \cite{kuruzov_stonyakin}, there is an analysis of the SESOP method with the required conditions on the quality of the solutions of the low-dimensional problems; these results are listed in Appendix \ref{appendix_sesop} for the reader's convenience. Using them, we are able to estimate the complexity of the SESOP method. Below, $D_k$ denotes the $n\times 3$ matrix with the columns $\mathbf{d}_k^0, \mathbf{d}_k^1, \mathbf{d}_k^2$, and $f_k(\tau) := f(x_k + D_k\tau)$ denotes the objective of the $k$-th auxiliary subproblem.

\begin{theorem}
\label{theorem_exact_tau}
Let the objective function $f$ be $L$-smooth and $\gamma$-quasar-convex with respect to $x^*$, and let the inexact gradient $g$ satisfy condition \eqref{g_cond_i} with $\delta_1\leq \frac{\varepsilon}{\frac{R}{\gamma}+10}$, where $R = \|x_0 - x^*\|$. Also, let us assume that on each iteration $k$ there is a ball $\mathcal{B}_{R_\tau}^k\subset \mathbb{R}^3$ of radius $R_\tau$ centered at zero such that $\tau_k \in \mathcal{B}_{R_\tau}^k$ and the subproblems are convex on these balls. If we can use the exact gradient of the functions $f_k$, then to approach the quality $\varepsilon$ for the initial problem by the SESOP method, one requires not more than 
    $$N = \mathcal{O}\left(\sqrt{\frac{LR^2}{\gamma^2 \varepsilon}}\right)$$
    inexact gradient calculations with respect to $x$ and not more than
    $$M=\mathcal{O}\left(N \ln \frac{L B C_N}{\varepsilon^4 }\right)$$
    exact gradient calculations with respect to $\tau$, where
    $$B = \max_{k=\overline{1, N}} \left(\max_{\tau\in\mathcal{B}_{R_\tau}^k} f_k(\tau) - \min_{\tau\in\mathcal{B}_{R_\tau}^k} f_k(\tau)\right),$$
    $$C_N = 1 + \left(\sqrt{\max_{k=\overline{1, N}}(\|D_k\| \|\tau_k\|)} + \sqrt{\max_{k=\overline{1, N}}\|\mathbf{d}^1_{k-1}\|}\right)^2+\max_{k=\overline{1, N}}\|\mathbf{d}_k^2\|.$$
\end{theorem}
\begin{proof}
    See detailed proof in Appendix~\ref{proof_theorem_exact_tau}.
\end{proof}

\begin{remark}
All the quantities entering $C_N$ are available during the run of the algorithm: the matrix $D_k$ is formed (in step 3 of Algorithm \ref{alg:sesop}) before the $k$-th subproblem is solved, and $\|\tau_k\| \leq R_\tau$ by assumption, so $\|D_k\|\|\tau_k\| \leq \|D_k\| R_\tau$. Thus, the required accuracy $\delta_4$ of the $k$-th subproblem (see Lemma \ref{delta_4_cond}) can be computed before this subproblem is solved. Moreover, $C_N$ is bounded whenever the trajectory $\{x_k\}$ of the method is bounded, since the columns of $D_k$ are expressed in terms of the inexact gradients at the trajectory points and the differences $x_k - x_0$.
\end{remark}

\begin{remark}
Note that the SESOP method in such an implementation requires $O\left(\sqrt{\frac{LR^2}{\gamma^2\varepsilon}}\right)$ inexact gradient calculations with respect to $x$ and $O\left(\sqrt{\frac{LR^2}{\gamma^2\varepsilon}}\ln \frac{LBC_N}{\varepsilon}\right)$ exact gradient calculations with respect to $\tau$.
\end{remark}

\begin{remark}
The main theoretical advantage of SESOP with an inexact gradient is that, in contrast to earlier results on accelerated methods with an inexact gradient (see, e.g., \cite{DevolderThesis, Vasin2021}), the error term in the convergence estimate (see Theorem \ref{SESOP_theorem_full}) does not contain an additive part proportional to $\max_k R_k$, where $R_k = \|x_k - x^*\|$ is the distance from the $k$-th trajectory point to the exact solution. This property is achieved through solving the additional low-dimensional subproblems. Nevertheless, it leads to requirements for a high accuracy of the solutions of the auxiliary problems.
\end{remark}

Further, let us consider the case of an inexact gradient in the internal problems. Namely, we assume that the $k$-th subproblem is solved by the ellipsoid method that uses the inexact low-dimensional gradient
$$g_k(\tau) = D_k^\top g(x_k + D_k\tau),$$
for which $\|g_k(\tau) - \nabla f_k(\tau)\| \leq \|D_k\|_2\delta_1$ (see Appendix \ref{proof_ellipsoids_sesop_final_theorem}). In this case, the quality of the subproblem solution cannot be better than the inexactness of the gradient (see \cite{Kuruzov_Gladin}). Uniting the results from \cite{Kuruzov_Gladin} for the ellipsoid method and the analysis of the previous theorem, we obtain the following estimate.

\begin{theorem}\label{ellipsoids_sesop_final_theorem}
Let the assumptions of Theorem \ref{theorem_exact_tau} hold with the following change: instead of the exact gradients of $f_k$, the subproblems are solved by the ellipsoid method with the inexact low-dimensional gradients $g_k(\tau) = D_k^\top g(x_k + D_k\tau)$, and the inexact gradient $g$ satisfies condition \eqref{g_cond_i} with \begin{equation}
\label{strong_delta_1}
    \delta_1\leq \min\left\{\frac{\varepsilon}{\frac{R}{\gamma}+10}, \frac{\varepsilon^4}{320000\, A_N L}\right\}
\end{equation} where $A_N=C_N \max_{k=\overline{1, N}} \|D_k\|_2$ for $C_N$ defined as in Theorem \ref{theorem_exact_tau}.
Then to approach the quality $\varepsilon$ for the initial problem by the SESOP method, one requires not more than
    $$N = \mathcal{O}\left(\sqrt{\frac{LR^2}{\gamma^2 \varepsilon}}\right)$$
    inexact gradient calculations with respect to $x$ and not more than
    $$M=\mathcal{O}\left(N \ln \frac{L B C_N}{\varepsilon^4 }\right)$$
    inexact gradient calculations with respect to $\tau.$
\end{theorem}
\begin{proof}
    See detailed proof in Appendix~\ref{proof_ellipsoids_sesop_final_theorem}.
\end{proof}
We can see that the number of gradient calculations is almost the same as in the previous theorem (the case of an exact low-dimensional gradient), but now we have the significantly stronger condition \eqref{strong_delta_1}. This condition guarantees inequalities \eqref{g_orth_i_eps} and \eqref{g_orth_ii_eps} and allows one to approach the quality \eqref{delta_4_eps} in the subproblem. But in this case, we can see that the gradient inexactness should be not more than $O(\varepsilon^4).$

Note that the advantages of the SESOP method lead to the requirements of an extra low inexactness both for the gradient and for the solutions of the auxiliary problems. Nevertheless, the required inexactness can be easily controlled through the values $\mathbf{d}_k^j$, $j = 0, 1, 2$, during the work of the algorithm (see the remark after Theorem \ref{theorem_exact_tau}).

\subsection{Multidimensional Dichotomy}

The ellipsoid method can be applied to the subproblems of the CG method too. Nevertheless, in \cite{Kuruzov_Gladin} it was shown that there is another effective method for the two-dimensional subproblem. It is a generalization of the one-dimensional dichotomy. Despite a slightly worse convergence rate in comparison with the ellipsoid method, it demonstrates better performance in practice. So, we provide the result for the CG method with the two-dimensional dichotomy in the following theorem. Similarly to the proof of Theorem~\ref{ellipsoids_sesop_final_theorem}, we are able to join the result of \cite{Kuruzov_Gladin} for the dichotomy method and Theorem \ref{thm-main} to obtain the following result.
\begin{theorem}\label{thm-main_halve_cube}
Let the assumptions of Theorem \ref{thm-main} hold with the following change: all the subproblems are convex and each point $\hat{x}_k$ is the output of the two-dimensional dichotomy algorithm (see \cite{Kuruzov_Gladin}) with the inexact low-dimensional gradient. Besides, there is $R_x$ such that $\hat{x}_k - x_k \in B_{R_x}$ for all $k.$ The number of steps $M$ of the dichotomy per subproblem is given by:
$$M=\left\lceil16 \left(\ln \frac{CR_x}{\varepsilon^4}\right)^2 \right\rceil,$$
where $C$ is a constant defined by the trajectory of the method similarly to $C_N$ in Theorem \ref{theorem_exact_tau}.

Let one of the following alternatives hold:
\begin{enumerate}
    \item The stopping criterion is never satisfied and the restarted Nemirovski's Conjugate Gradient method (Algorithm \ref{alg:CG_restarts}) makes $$K = \left\lceil\frac{2}{1-\alpha}\log\frac{\varepsilon_0}{\varepsilon}\right\rceil$$
restarts with
$$T=\left\lceil\frac{4}{\gamma} \sqrt{\frac{L}{\mu}}\frac{\sqrt{1+\alpha}}{1-\alpha}\right\rceil+1$$
iterations per each restart, where $\varepsilon=\frac{196}{\gamma^2 \mu} \delta_1^2$ and $\alpha\in\left[\frac{1}{\sqrt 2},1\right)$;
\item At some inner iteration, the stopping criterion $\|g(\hat{x}_N)\|\leq \frac{8}{\gamma}\delta_1$ is satisfied for the first time at a computed point $\hat{x}_N$.
\end{enumerate}
Then for the output point $\widehat{x}$ ($\widehat{x} = x_{K}$ in the first case and $\widehat{x} = \hat{x}_{N}$ in the second case), the following inequality holds:
$$
f(\widehat{x})-f^{*}\leqslant\frac{196\delta_1^{2}}{\gamma^2 \mu}.
$$
As a result, the algorithm requires not more than $KT$ calculations of the inexact gradient with respect to $x$ and $MKT=O(\ln^3 (1/\varepsilon))$ calculations of the low-dimensional inexact gradient.
\end{theorem}

\begin{proof}
The alternatives and the final accuracy bound are established exactly as in the proof of Theorem \ref{thm-main} (using the implication PL $\Rightarrow$ quadratic growth from \cite{Karimi} and Theorem \ref{CG_theorem_restarts}), since the two-dimensional dichotomy with the number of steps $M$ from the statement provides the accuracy of the subproblem solution sufficient for the (approximate) orthogonality conditions used in the proofs; the estimate for $M$ per subproblem follows from the complexity result for the two-dimensional dichotomy with an inexact gradient from \cite{Kuruzov_Gladin}, applied with the required subproblem accuracy $\delta_4 = O\left(\frac{\varepsilon^4}{LC}\right)$ obtained similarly to Lemma \ref{delta_4_cond}.
\end{proof}

We can see that the asymptotic complexity of the CG method with the dichotomy is proportional to the cube of $\ln (1/\varepsilon).$ On the one hand, it is clear that the ellipsoid method gives a better asymptotic complexity. On the other hand, \cite{Kuruzov_Gladin} demonstrates that the dichotomy method outperforms the ellipsoid method for two-dimensional problems. So, let us consider the performance of the considered approaches in practice.

\section{Numerical Experiments}\label{section:numer_exp}

In this section, we present some numerical results. We compare the SESOP method with low-dimensional methods for the auxiliary subproblems (see Algorithms \ref{alg:sesop} and \ref{alg:ellispoids}), the CG method with restarts (see Algorithms \ref{alg:CG_restarts}, \ref{alg:ellispoids} and \cite{Kuruzov_Gladin}), the Similar Triangle Method (STM, see \cite{stm}) and Gradient Descent (GD) as baselines.

{\bf Logistic regression (a convex baseline for the time comparison).} First, as a convex baseline, we consider the problem of regularized logistic regression:
\begin{equation}\label{log_reg}
    f(x)=({1}/{m})\sum_{j=1}^m\log(1+\exp(-y_{j}\langle a_{j}, x\rangle)) + \mu_{reg} \|x\|^2
\end{equation}
on synthetic data with dimensions $n=100$, $m=200$. We include this convex example specifically because here the objective is convex, and therefore all the low-dimensional auxiliary subproblems (the three-dimensional subproblems of SESOP and the two-dimensional subproblems of the CG method) are guaranteed to be convex. This is exactly the setting in which the complexity guarantees of Section \ref{section:subproblems} apply and the low-dimensional methods (the ellipsoid method and the multidimensional dichotomy) are known to converge correctly; the purpose of this experiment is thus to compare the running time of the considered approaches in a regime that is fully covered by the theory. The data were generated as follows: the features $a_j \in \mathbb{R}^n$ were drawn from the standard normal distribution $\mathcal{N}(0, I_n)$, and the labels were set as $y_j = \mathrm{sign}(\langle a_j, \bar{x}\rangle)$ for a random vector $\bar{x}\sim\mathcal{N}(0, I_n)$. We used synthetic data (rather than, e.g., the LibSVM datasets) to be able to control all the parameters of the problem and to report the results averaged over several independently generated datasets. We set the regularization parameter $\mu_{reg} = 10^{-3}$, so the problem is $\mu$-strongly convex with $\mu = 2\mu_{reg}$; the Lipschitz constant was computed analytically as $L = \frac{1}{4m}\lambda_{\max}(A^\top A) + 2\mu_{reg}$, where $A$ is the matrix with the rows $a_j^\top$. The parameter of the restarts for the CG method was chosen according to Theorem \ref{CG_theorem_restarts}. The inexact gradient was modeled as $g(x) = \nabla f(x) + \delta_1 e$, where $e$ is a fixed random unit vector, i.e. we used a deterministic additive perturbation of the norm exactly $\delta_1$; note that this perturbation is not stochastic, in accordance with the considered inexactness model \eqref{InexactGrad}. It was found that all the methods approach close accuracies in our problem settings, as expected in this convex regime. So, we demonstrate the time comparison (see Table \ref{table:res}); the presented time (averaged over 10 generated datasets) is the time required to approach the quality $10\delta_1^2/\mu$ in the function value.

In the remaining experiments we move beyond this convex baseline to a non-convex quasar-convex problem, consider both an ill-conditioned and a well-conditioned regime, and add GD and STM to the comparison.

\begin{table}[ht]
    \centering
    \begin{tabular}{|c|c|c|c|c|}
         \hline
         $\delta_1$& SESOP & CG+Ellipsoids & CG+Dichotomy & STM \\
         \hline
         $10^{-3}$ & 1 & 1.4 & 0.9& 1.7\\
         $10^{-5}$ & 10.1 & 15.3& 9.5& 13.8\\
         $10^{-7}$ & 35.3 & 60.9 & 36.8& 42.1\\
         \hline
    \end{tabular}
    \caption{Time comparison (s) for problem \eqref{log_reg}}
    \label{table:res}
\end{table}

We can see that the multidimensional dichotomy works better than the ellipsoid method for Nemirovski's Conjugate Gradient method in all cases. At the same time, the methods based on subspace optimization outperform the STM method. The results for the CG and SESOP methods are close enough.

\begin{figure}[ht]
    \centering
    \includegraphics[width=0.9\textwidth]{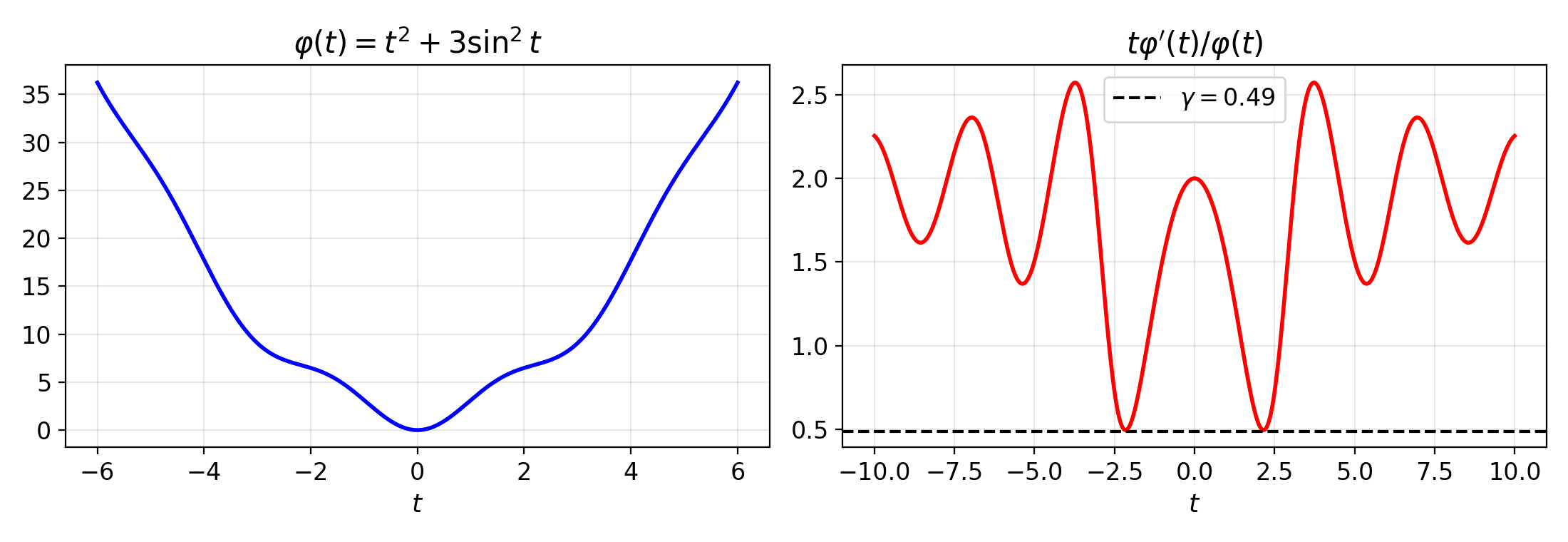}
    \caption{The one-dimensional profile $\varphi(t) = t^2 + 3\sin^2 t$ of the non-convex problem \eqref{a_composed_problem} (left) and the ratio $t\varphi'(t)/\varphi(t)$ (right), which is separated from zero by $\gamma = 0.49$ (Appendix \ref{appendix_quasar_example}).}
    \label{fig:quasar_check}
\end{figure}

{\bf A non-convex test problem.} To validate the methods on the intended function class beyond the convex baseline, in the remaining experiments we consider the non-convex problem
\begin{equation}\label{a_composed_problem}
f(x) \;=\; \sum_{i=1}^{n}\Big( (a_i^\top x + b_i)^2 + 3\sin^2(a_i^\top x + b_i)\Big) \;=\; \|Ax+b\|^2 + 3\sum_{i=1}^n \sin^2\big((Ax+b)_i\big),
\end{equation}
where $a_i^\top$ are the rows of an invertible matrix $A \in \mathbb{R}^{n\times n}$ and $b \in \mathbb{R}^n$. The unique stationary point is $x^* = A^{-1}(-b)$ with $f^* = 0$. Each term is built from the same scalar profile $\varphi(t) = t^2 + 3\sin^2 t$ used in the example \eqref{example_pl_quasar}; the composition with $A$ makes the Hessian of $f$ indefinite-directional, so $f$ is non-convex. As shown in Appendix \ref{appendix_quasar_example}, $f$ is $\gamma$-quasar-convex with respect to $x^*$ with $\gamma \geq 0.49$, and this constant does not depend on $A$: with $y = Ax + b$ one has $\langle \nabla f(x), x - x^*\rangle = h(y)^\top y = \sum_i y_i \varphi'(y_i)$ (where $h_i(y) = \varphi'(y_i)$), so the quasar-convexity of $f$ reduces to the {\it same} one-dimensional bound $t\varphi'(t) \geq \gamma\varphi(t)$ as for the separable example, whose exact constant is $\approx 0.496$ (see Figure \ref{fig:quasar_check}). Moreover, again independently of $A$, the same reduction shows that $f$ inherits the PL-condition of the scalar profile: $f$ satisfies the PL-condition with $\mu_{PL} = \sigma_{\min}^2(A)/32$ and the quadratic growth condition with $\mu = 2\sigma_{\min}^2(A)$, while $f$ is $L$-smooth with $L = 8\sigma_{\max}^2(A)$; hence the condition number of the problem is $\kappa = L/\mu = 4\big(\sigma_{\max}(A)/\sigma_{\min}(A)\big)^2$ and is fully controlled by the spectrum of $A$. In the experiments below, $A = U\,\mathrm{diag}(s)\,V^\top$ with random orthogonal $U, V$ and the singular values $s_i$ spread uniformly over $[1, \sigma_{\max}]$, and $b$ is a standard Gaussian vector.

\begin{figure}[ht]
    \centering
    \includegraphics[width=0.99\textwidth]{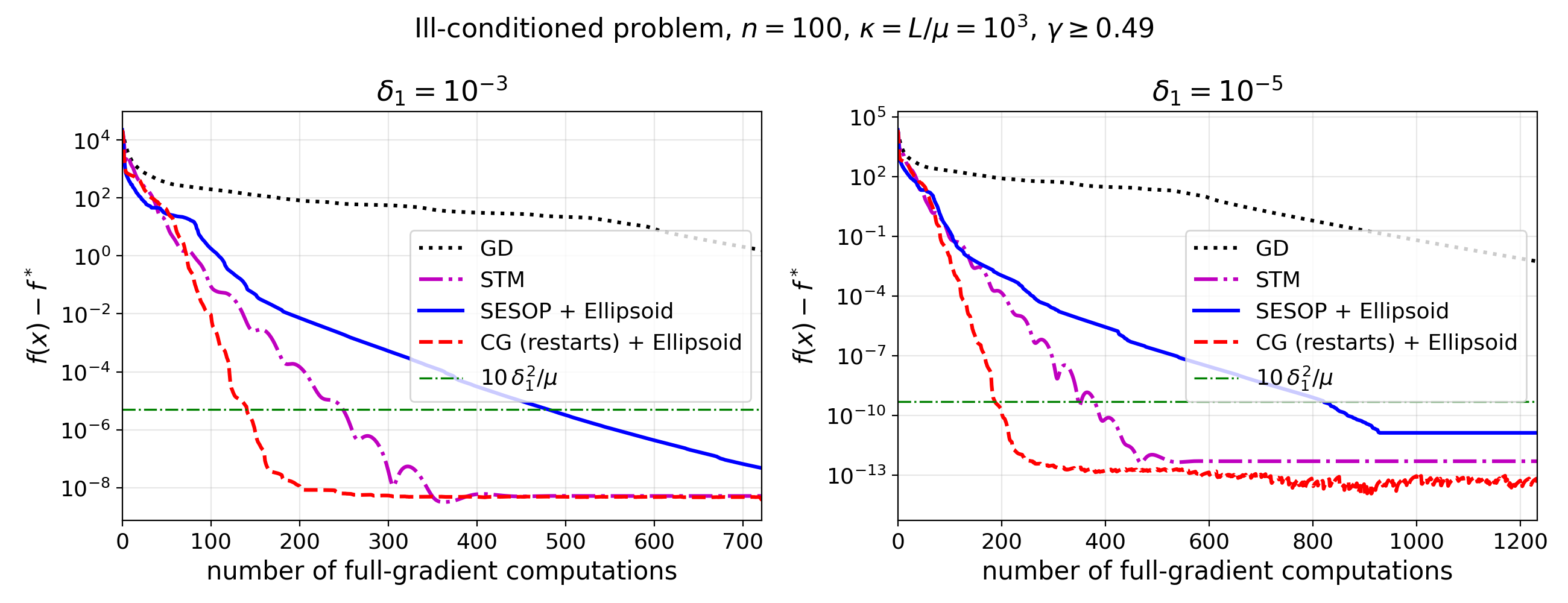}
    \caption{Convergence of GD, STM, SESOP and the restarted CG method on the ill-conditioned non-convex problem \eqref{a_composed_problem} ($n = 100$, $\kappa = 10^3$, $\gamma \geq 0.49$) with respect to the number of full-gradient computations, for the noise levels $\delta_1 = 10^{-3}$ (left) and $\delta_1 = 10^{-5}$ (right).}
    \label{fig:ill_cond}
\end{figure}

\begin{figure}[ht]
    \centering
    \includegraphics[width=0.99\textwidth]{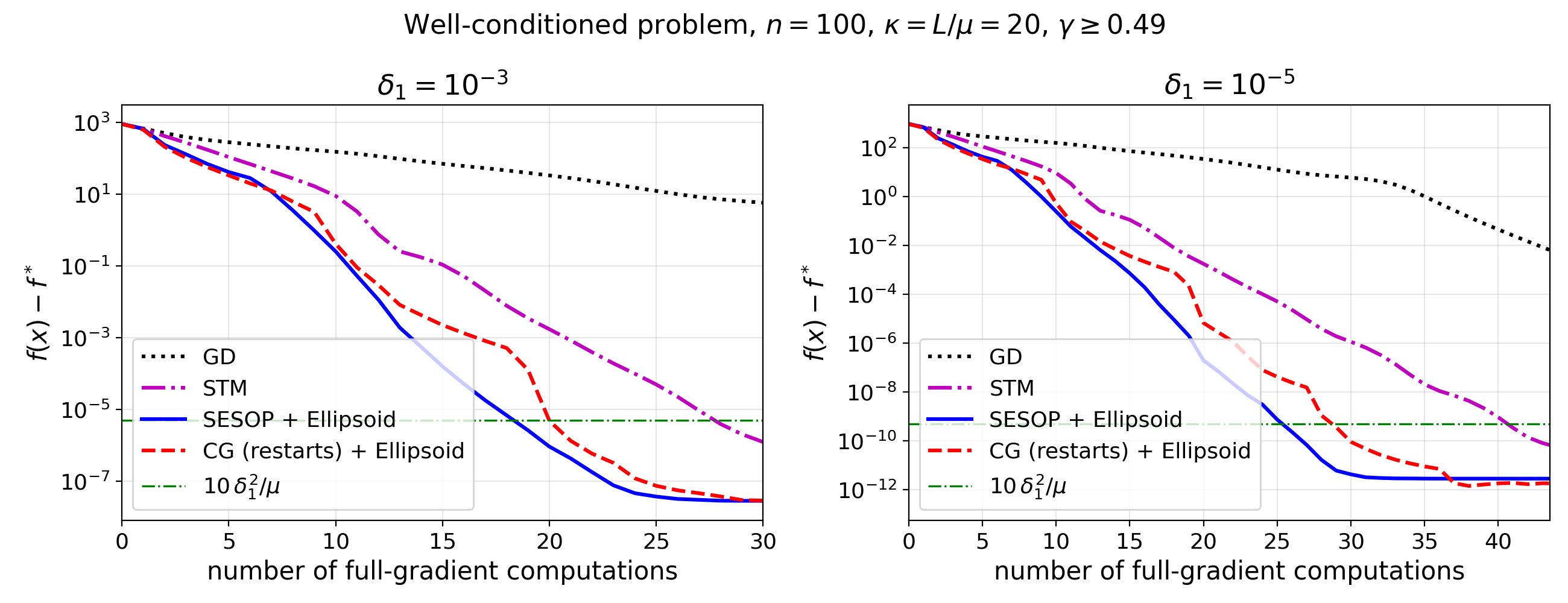}
    \caption{Convergence of GD, STM, SESOP and the restarted CG method on the well-conditioned non-convex problem \eqref{a_composed_problem} ($n = 100$, $\kappa = 20$, $\gamma \geq 0.49$) with respect to the number of full-gradient computations, for the noise levels $\delta_1 = 10^{-3}$ (left) and $\delta_1 = 10^{-5}$ (right).}
    \label{fig:well_cond}
\end{figure}

{\bf Measurement methodology.} Starting from this point, all the methods are compared with respect to the number of {\it full-dimensional} gradient computations, which for the problem \eqref{a_composed_problem} is the natural unit of cost: one full gradient $\nabla f(x) = A^\top(2y + \sin 2y)$, $y = Ax+b$, requires two matrix--vector products with $A$, i.e. $O(n^2)$ arithmetic operations. A key observation is that the subspace methods (SESOP and the restarted CG method) spend most of their oracle calls on the {\it low-dimensional} gradients of the auxiliary subproblems, and these are much cheaper than the full-dimensional ones: once the vectors $Ax+b$ and $AD_k$ are available, each low-dimensional gradient costs only $O(n)$ operations. Moreover, with a simple caching strategy the columns of $AD_k$ come essentially for free, because the directions used by both methods ($x_k - x_0$ and the weighted sums $q_k$ of past gradients) are linear combinations of vectors whose products with $A$ have already been computed at the previous iterations. As a result, one outer iteration of SESOP or of the restarted CG method costs the same two matrix--vector products as one iteration of GD or STM, i.e. exactly one full-gradient computation; the detailed operation count, together with a discussion of the (modest, two to four extra vectors) memory overhead, is given in Appendix \ref{appendix_cost}. This is why the horizontal axes of Figures \ref{fig:ill_cond} and \ref{fig:well_cond} report the number of full-gradient computations, which coincides with the number of outer iterations for all four methods. Let us also note that, in contrast with the logistic regression experiment, in the remaining experiments we use the ellipsoid method (rather than the multidimensional dichotomy) as the {\it single} low-dimensional solver for both SESOP and the restarted CG method: this keeps the comparison uniform --- both subspace methods are then built from exactly the same computational blocks and differ only in the choice of the subspaces --- and allows us to focus on the per-iteration convergence of the outer schemes themselves rather than on the specifics of a particular low-dimensional solver. The subproblems were solved with a fixed small number $N_{ell} = 10$ of ellipsoid iterations and an adaptively updated localization radius $R_\tau$; these practical heuristics are applied to both subspace methods on an equal footing, and Appendix \ref{appendix_nell} demonstrates that the restarted CG method is robust with respect to the choice of $N_{ell}$. The number of iterations per restart of the CG method was set to $T = 20$; the theoretical value $T = \lceil 2\sqrt{L/\mu}\rceil$ of Theorem \ref{CG_theorem_restarts} guarantees the linear rate, but in practice a smaller restart period turned out to be somewhat faster on this problem. The inexact gradient was modeled, as before, by a deterministic additive perturbation of norm exactly $\delta_1$ along a fixed random unit direction.

{\bf An ill-conditioned problem.} We first take $n = 100$ and $\sigma_{\max} = \sqrt{\kappa/4}$ with $\kappa = 10^3$. Figure \ref{fig:ill_cond} compares GD, STM, SESOP and the restarted CG method for $\delta_1 \in \{10^{-3}, 10^{-5}\}$; the horizontal axis is clipped at $1.5$ times the largest number of full-gradient computations needed by the two subspace methods to reach the threshold $10\delta_1^2/\mu$. The results agree well with the theory. The restarted CG method is the fastest of the four: it reaches the noise floor in about $140$--$190$ full-gradient computations, ahead of the accelerated method STM ($250$--$350$) and far ahead of SESOP ($480$--$820$) and GD ($2000$--$3100$). The advantage over GD (more than an order of magnitude) is the direct consequence of the linear rate $O\big(\frac{1}{\gamma}\sqrt{L/\mu}\log(1/\varepsilon)\big)$ of Theorem \ref{thm-main} versus the $O(\kappa\log(1/\varepsilon))$ rate of GD on an ill-conditioned problem; the advantage over SESOP reflects the linear-versus-sublinear comparison of the guarantees of Sections \ref{section:conj_grad} and \ref{section:subproblems}. All four methods converge to the same error floor of order $O(\delta_1^2/\mu)$, which again confirms that the gradient inexactness does not accumulate.

{\bf A well-conditioned counterpart.} For comparison, we repeat the experiment on the same problem \eqref{a_composed_problem} with a well-conditioned matrix, $\kappa = 20$, keeping $n = 100$ and the same noise levels. Figure \ref{fig:well_cond} shows that in this regime all the gaps shrink, as expected: both subspace methods reach the noise floor in $19$--$29$ full-gradient computations, on par with (and for SESOP slightly ahead of) STM ($28$--$41$), while GD remains the slowest ($58$--$77$). The two subspace methods are thus competitive in both regimes, and the restarted CG method is clearly preferable on ill-conditioned problems, in full agreement with its linear convergence guarantee.

\section*{Conclusion}

In this paper, we considered a generalization of the convexity condition that is known as quasar-convexity or weak quasi-convexity. We proposed a modification of Nemirovski's Conjugate Gradient method with a $\delta_1$-additive noise in the gradient \eqref{InexactGrad} for $\gamma$-quasar-convex functions satisfying the quadratic growth condition.

We estimated the computational complexity of solving the internal auxiliary problems in the SESOP and CG methods. For this, we used well-known low-dimensional optimization methods --- the ellipsoid method and a generalization of the dichotomy. We proved that these methods do not significantly increase the complexity in the case of convex subproblems. Besides, the resulting methods are still first-order methods.

Moreover, we provided numerical experiments, including experiments on a non-convex quasar-convex problem, which demonstrate the effectiveness of the approaches proposed in this paper. In particular, our experiments on a non-convex quasar-convex problem show that, when the cost is measured in full-dimensional gradient computations (which is the natural metric for the considered problem class, since the subspace methods mostly rely on much cheaper low-dimensional gradients; see Appendix \ref{appendix_cost}), the restarted CG method is the fastest of the compared methods on an ill-conditioned problem --- ahead of the accelerated method STM and more than an order of magnitude ahead of GD --- owing to its linear convergence rate, while on a well-conditioned problem both subspace methods remain competitive with STM; in all the considered regimes the methods converge to the same noise floor of order $O(\delta_1^2/\mu)$ without accumulating the gradient inexactness.

\bibliography{biblio_new}

\newpage
\appendix

\section{Quasar-convexity of the example \eqref{example_pl_quasar}}
\label{appendix_quasar_example}
Consider the one-dimensional function $\varphi(t) = t^2 + 3\sin^2 t$, so that the separable function \eqref{example_pl_quasar} is $f(x) = \sum_{i=1}^n \varphi(x_i)$. Since $f$ is separable and $x^* = 0$ with $f^* = 0$, the $\gamma$-quasar-convexity of $f$ with respect to $x^* = 0$ follows from the one-dimensional inequality
\begin{equation}\label{quasar_1d}
\varphi(x^*) \geq \varphi(t) + \frac{1}{\gamma}\varphi'(t)(x^* - t)\Big|_{x^*=0} = \varphi(t) - \frac{1}{\gamma}t\varphi'(t),
\end{equation}
i.e. from $t\varphi'(t) \geq \gamma\,\varphi(t)$ for all $t$ (indeed, summing \eqref{quasar_1d} over the coordinates gives the definition of quasar-convexity for $f$; recall that $\varphi(t)\geq 0$). Here $\varphi'(t) = 2t + 3\sin 2t$. We prove that
\begin{equation}\label{g_positive}
G(t) := t\varphi'(t) - \gamma\,\varphi(t) = (2-\gamma)t^2 + 3t\sin 2t - 3\gamma\sin^2 t \geq 0\qquad\forall t\in\mathbb{R}
\end{equation}
for $\gamma = 0.49$. The function $G$ is even, so it suffices to consider $t > 0$. Substituting $\gamma = 0.49$, we get $G(t) = 1.51\,t^2 + 3t\sin 2t - 1.47\sin^2 t$. We split the half-line into three regions.

\medskip
\noindent\textit{Region 1: $0 < t \leq 1.22$.} Using the elementary inequalities $\sin 2t \geq 2t - \frac{(2t)^3}{6}$ and $\sin^2 t \leq t^2$, we obtain $G(t) \geq t^2\left(6.04 - 4t^2\right)$, which is strictly positive for $t \leq 1.22 < \sqrt{6.04/4}$.

\medskip
\noindent\textit{Region 2: $1.22 \leq t \leq 2.40$.} On this compact interval we use a Lipschitz argument. Since $G'(t) = 3.02\,t + 1.53\sin 2t + 6t\cos 2t$, we have $|G'(t)| \leq 23.18$ for $t\in[1.22, 2.40]$. Evaluating $G$ on a uniform grid with step $0.002$ gives $G(t_j) \geq 0.041$ at every node, and since any point of the interval is within $0.001$ of a node, $G(t) > 0$ on the whole interval.

\medskip
\noindent\textit{Region 3: $t \geq 2.40$.} Using $\sin 2t \geq -1$ and $\sin^2 t \leq 1$, we get $G(t) \geq 1.51\,t^2 - 3t - 1.47 =: H(t)$. The parabola $H$ is increasing on $[2.40, \infty)$ and its larger root is $\approx 2.394 < 2.40$, so $H(t) \geq H(2.40) > 0$ there.

\medskip
Combining the three regions, $G(t) \geq 0$ for all $t \in \mathbb{R}$, so $\varphi$ (and therefore $f$) is $0.49$-quasar-convex with respect to $x^* = 0$. The exact value of $\inf_{t\neq 0} t\varphi'(t)/\varphi(t)$ is $\approx 0.496$ (attained near $t \approx \pm 2.15$), so the bound $\gamma = 0.49$ is close to tight.

\medskip
Consider now the non-convex problem \eqref{a_composed_problem} used in the numerical experiments of Section \ref{section:numer_exp},
$$f(x) = \|Ax+b\|^2 + 3\sum_{i=1}^n \sin^2\big((Ax+b)_i\big) = \sum_{i=1}^n \varphi\big((Ax+b)_i\big),$$
with an invertible $A$ and the {\it same} scalar profile $\varphi(t) = t^2 + 3\sin^2 t$ as above; let $x^* = A^{-1}(-b)$, $f^* = 0$. Writing $y = Ax + b$, we have $\nabla f(x) = A^\top h(y)$ with $h_i(y) = \varphi'(y_i)$, and since $x - x^* = A^{-1}y$,
$$\langle \nabla f(x),\, x - x^*\rangle = h(y)^\top A A^{-1} y = h(y)^\top y = \sum_{i=1}^n y_i\,\varphi'(y_i).$$
Thus the quasar-convexity of $f$ reduces to exactly the same one-dimensional inequality $t\varphi'(t) \geq \gamma\varphi(t)$ that was proved above; in particular, {\it the constant does not depend on $A$}, and $f$ is $\gamma$-quasar-convex with respect to $x^*$ with the same $\gamma = 0.49$ (exact value $\approx 0.496$; see Figure \ref{fig:quasar_check}). There is therefore no need to repeat the three-region argument.

The same separable reduction transfers the other properties of $\varphi$ to $f$, again independently of $A$:
\begin{itemize}
\item[(i)] {\it PL-condition.} The scalar profile satisfies $\varphi'(t)^2 \geq 2\mu_0\,\varphi(t)$ with $\mu_0 = 1/32$ (this is the one-dimensional PL-constant of $x^2 + 3\sin^2 x$ used in \eqref{example_pl_quasar}, see \cite{Karimi}). Hence $\|\nabla f(x)\|^2 = h(y)^\top A A^\top h(y) \geq \sigma_{\min}^2(A)\,\|h(y)\|^2 = \sigma_{\min}^2(A)\sum_i \varphi'(y_i)^2 \geq 2\sigma_{\min}^2(A)\mu_0 \sum_i \varphi(y_i)$, so $f$ satisfies the PL-condition with constant $\mu_{PL} = \sigma_{\min}^2(A)/32$.
\item[(ii)] {\it Quadratic growth.} Since $\varphi(t) \geq t^2$, $f(x) - f^* \geq \|y\|^2 = \|A(x-x^*)\|^2 \geq \sigma_{\min}^2(A)\,\|x-x^*\|^2$, so the quadratic growth condition \eqref{QC_cond} holds with $\mu = 2\sigma_{\min}^2(A)$.
\item[(iii)] {\it Smoothness.} As $\nabla^2 f(x) = A^\top\,\mathrm{diag}(\varphi''(y_i))\,A$ with $\varphi''(t) = 2 + 6\cos 2t \in [-4, 8]$, we have $\nabla^2 f(x) \preceq 8\,A^\top A$, i.e. $f$ is $L$-smooth with $L = 8\sigma_{\max}^2(A)$.
\end{itemize}
Consequently the condition number of the problem is $\kappa = L/\mu = 4\big(\sigma_{\max}(A)/\sigma_{\min}(A)\big)^2$, controlled entirely by the spectrum of $A$.

\section{Cost of one iteration of the subspace methods}
\label{appendix_cost}
In this appendix we justify the measurement methodology of Section \ref{section:numer_exp}: for the problem \eqref{a_composed_problem}, one outer iteration of SESOP or of the restarted CG method costs, up to $O(n)$ terms, the same amount of arithmetic as one iteration of GD or STM, namely one full-gradient computation, i.e. two matrix--vector products with the matrix $A$ (each of cost $O(n^2)$).

{\bf Full gradient.} With $y = Ax + b$ we have $\nabla f(x) = A^\top\big(2y + \sin 2y\big)$, so one full gradient requires one product with $A$ (to form $y$) and one product with $A^\top$: two matrix--vector products in total. This is the per-iteration cost of GD and STM.

{\bf Low-dimensional gradients are cheap.} The objective of the auxiliary subproblem is $f_k(\tau) = f(x_k + D_k\tau)$ with $y(\tau) = (Ax_k + b) + (AD_k)\tau$. Once the vector $Ax_k + b$ and the (two or three columns of the) matrix $AD_k$ are available, each evaluation of the low-dimensional gradient $\nabla f_k(\tau) = (AD_k)^\top\big(2y(\tau) + \sin 2y(\tau)\big)$ costs only $O(n)$ operations. Hence all the $N_{ell}$ inner iterations of the ellipsoid method together cost $O(N_{ell}\, n) = o(n^2)$ and are negligible compared to a single matrix--vector product.

{\bf The columns of $AD_k$ come for free.} The directions used by the two subspace methods are $x_k - x_0$ and the aggregated gradient sums $q_k$ (plus, for SESOP, the current gradient $g(x_k)$). All of them are linear combinations of vectors whose products with $A$ are already known:
\begin{itemize}
\item $A(x_k - x_0) = (Ax_k + b) - (Ax_0 + b)$, where $Ax_k + b$ is maintained along the trajectory (see below) and $Ax_0 + b$ is computed once at the start;
\item $q_{k+1} = q_k + w_k\, g(\hat{x}_k)$, so $Aq_{k+1} = Aq_k + w_k\, A g(\hat{x}_k)$, i.e. the cached vector $Aq_k$ is updated by the single product $A g(\hat{x}_k)$;
\item after the subspace step $x_{k+1} = \hat{x}_k - \frac{1}{2L} g(\hat{x}_k)$ (or $x_{k+1} = x_k + D_k \tau_k$ for SESOP), the vector $Ax_{k+1} + b$ is obtained from the cached quantities: $A\hat{x}_k + b = (Ax_k + b) + (AD_k)\tau_k$ is already available from the subproblem, and $A x_{k+1} + b = (A\hat{x}_k + b) - \frac{1}{2L} A g(\hat{x}_k)$ reuses the same product $A g(\hat{x}_k)$.
\end{itemize}
Therefore one outer iteration of either subspace method requires exactly two matrix--vector products beyond the cached data: the product with $A^\top$ inside the full gradient $g(\hat{x}_k)$, and the product $A g(\hat{x}_k)$ used to refresh the caches. This is the same count as for one full gradient, so measuring all four methods in the number of full-gradient computations (equivalently, outer iterations) is a fair comparison.

{\bf Memory.} The caching above requires storing a small constant number of additional $n$-dimensional vectors: $Ax_k + b$, $Ax_0 + b$, $Aq_k$ and (temporarily) $Ag(\hat{x}_k)$ --- no more than two to four extra vectors overall, and this bookkeeping can be organized even more economically by overwriting the temporaries in place. The memory overhead is thus $O(n)$ and negligible compared to the $O(n^2)$ storage of the matrix $A$ itself.

\section{Robustness of the CG method with respect to $N_{ell}$}
\label{appendix_nell}
\begin{figure}[ht]
    \centering
    \includegraphics[width=0.99\textwidth]{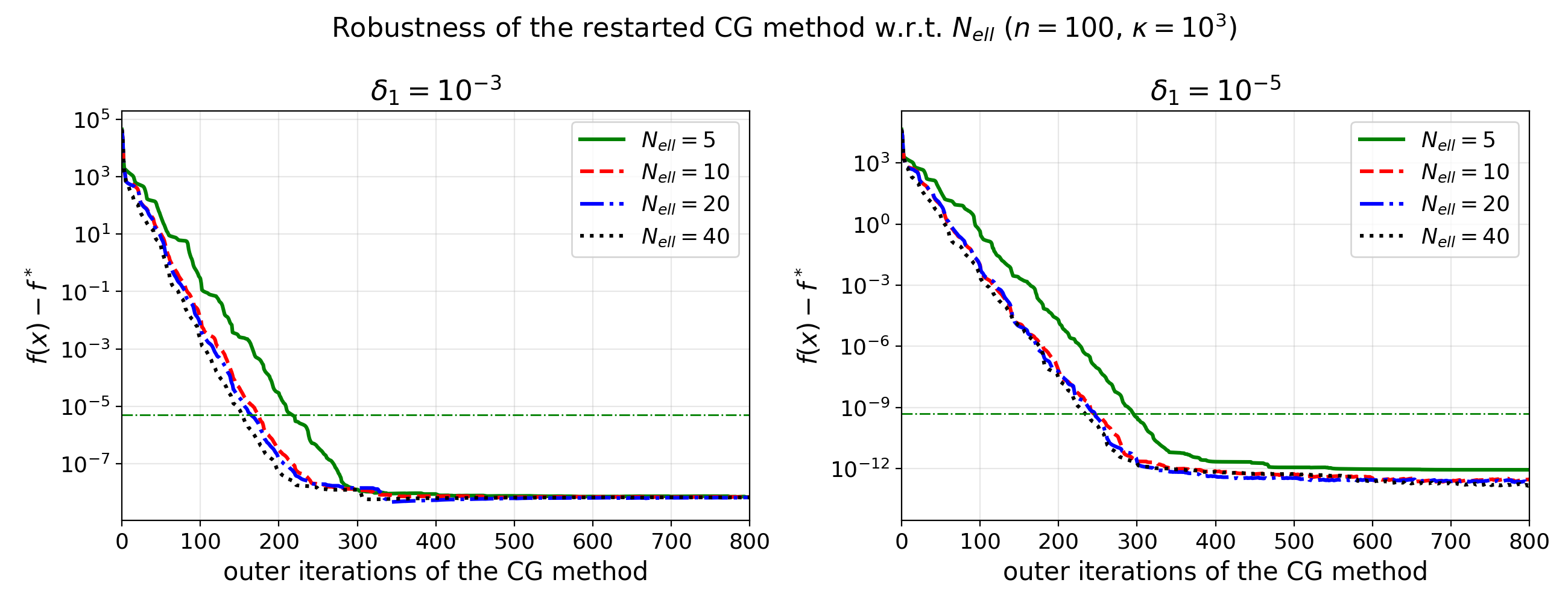}
    \caption{Convergence of the restarted CG method on the ill-conditioned problem \eqref{a_composed_problem} ($n = 100$, $\kappa = 10^3$) for different numbers $N_{ell}$ of inner ellipsoid iterations, with respect to the number of outer iterations, for $\delta_1 = 10^{-3}$ (left) and $\delta_1 = 10^{-5}$ (right).}
    \label{fig:nell}
\end{figure}

The experiments of Section \ref{section:numer_exp} use a fixed small number $N_{ell} = 10$ of inner iterations of the ellipsoid method per auxiliary subproblem. Figure \ref{fig:nell} shows the convergence of the restarted CG method on the ill-conditioned problem of Section \ref{section:numer_exp} ($n = 100$, $\kappa = 10^3$) for $N_{ell} \in \{5, 10, 20, 40\}$, with respect to the number of outer iterations. The trajectories are close to each other, and the number of outer iterations needed to reach the noise floor varies only mildly (for instance, for $\delta_1 = 10^{-3}$: about $198$ outer iterations for $N_{ell} = 5$, $141$ for $N_{ell} = 10$, and $122$ for both $N_{ell} = 20$ and $N_{ell} = 40$). At the same time, the total number of {\it low-dimensional} gradient calls grows linearly with $N_{ell}$, so larger values of $N_{ell}$ buy a marginal reduction of the outer iteration count at a proportionally growing (though cheap, see Appendix \ref{appendix_cost}) inner cost. We conclude that the restarted CG method is robust with respect to the accuracy of the inner solver, and a small fixed $N_{ell}$ is a reasonable practical choice.

\section{Proof of Lemma \ref{CG_lemma}}
\label{appendix_lemma_proof}
Note that, according to Algorithm \ref{alg:CG}, $q_{k+1} = q_{k} + g(\hat{x}_{k})$ for $k\geq 0$. So, we have the following expression for $\|q_{k+1}\|^2$:
$$\|q_{k+1}\|^2 = \|q_{k}\|^2 + \|g(\hat{x}_{k})\|^2 + 2\langle g(\hat{x}_{k}), q_{k}\rangle.$$
Since the auxiliary problem \eqref{eqsubproblem_cg} is solved exactly at each iteration, the point $\hat{x}_k$ is a minimizer of $f$ over the affine subspace $\mathcal{X}_k = x_0 + \mathrm{Lin}(x_k - x_0, q_k)$, and therefore the gradient $\nabla f(\hat{x}_k)$ is orthogonal to the directions of this subspace; in particular,
$$\nabla f (\hat{x}_k)\perp q_{k}.$$
At the same time, we have $\|\nabla f (\hat{x}_k)-g(\hat{x}_k)\|\leq \delta_1$ by \eqref{g_cond_i}. Hence,
$$|\langle g(\hat{x}_k), q_k\rangle| = |\langle g(\hat{x}_k) - \nabla f(\hat{x}_k), q_k\rangle| \leq \delta_1\|q_k\|.$$
So, summing the expression for $\|q_{k+1}\|^2$ over $k = 0, \ldots, T-1$, we obtain the upper bound
\begin{equation}
\label{upper_bound}
\|q_T\|^2 \leq \sum\limits_{k=0}^{T-1} \|g(\hat{x}_k)\|^2 + 2\delta_1 \sum\limits_{k=0}^{T-1} \|q_k\|,
\end{equation}
and, for each $k$, the lower bound
\begin{equation}
\label{lower_bound}
\|q_{k+1}\|^2 \geq \|q_{k}\|^2 - 2\delta_1\|q_{k}\|.
\end{equation}
From the lower bound \eqref{lower_bound}, we have that 
$$\|q_{k}\|\leq \delta_1 + \sqrt{\delta_1^2 + \|q_{k+1}\|^2}\leq 2\delta_1 + \|q_{k+1}\|.$$
Applying the inequality above repeatedly, we obtain an estimate of $\|q_j\|$ for all $j\leq T$:
\begin{equation}
\label{lower_bound_j}
\|q_j\|\leq 2\delta_1 (T-j) + \|q_T\|.
\end{equation}
Using inequalities \eqref{upper_bound} and \eqref{lower_bound_j}, we have the following estimate:
\begin{equation}
\label{upper_bound_new}
\|q_T\|^2 \leq \sum\limits_{k=0}^{T-1} \|g(\hat{x}_k)\|^2 + 4\delta_1^2 \sum\limits_{j=0}^{T-1} (T-j) + 2\delta_1 T \|q_T\|.
\end{equation}
Since $\sum_{j=0}^{T-1}(T-j) = \frac{T(T+1)}{2} \leq T^2$, from \eqref{upper_bound_new} we obtain the following quadratic inequality on $\|q_T\|$:
\begin{equation}
\label{upper_bound_new_1}
\|q_T\|^2 \leq \sum\limits_{k=0}^{T-1} \|g(\hat{x}_k)\|^2 + 4\delta_1^2 T^2 + 2\delta_1 T \|q_T\|.
\end{equation}
Solving the quadratic inequality \eqref{upper_bound_new_1} with respect to $\|q_T\|$ and using $\sqrt{a+b}\leq \sqrt a + \sqrt b$, we obtain
$$\|q_T\|\leq \delta_1 T + \sqrt{\delta_1^2T^2 + 4\delta_1^2T^2 + \sum\limits_{k=0}^{T-1} \|g(\hat{x}_k)\|^2} \leq (1+\sqrt{5})\delta_1 T + \left(\sum\limits_{k=0}^{T-1} \|g(\hat{x}_k)\|^2\right)^{\frac{1}{2}},$$
which gives estimate \eqref{q_t_est} since $1 + \sqrt 5 \leq 4$.

\section{Some Results for SESOP method}
\label{appendix_sesop}
Let us provide here some results from \cite{kuruzov_stonyakin} for the SESOP method that will be useful for us in the further proofs. Let us start from the basic inequality for the SESOP method. Since the direction $\mathbf{d}_k^0 = g(\mathbf{x}_k)$ belongs to the subspace over which the minimization is performed in step 4 of Algorithm \ref{alg:sesop}, we have
\begin{equation}
    \label{k_min_ineq}
    f(\mathbf{x}_{k+1}) = \min_{\mathbf{s}\in\mathbb{R}^3} f\left(\mathbf{x}_k + \sum\limits_{i=0}^2 s_i \mathbf{d}_k^i\right)\leq f\left(\mathbf{x}_k + s_0 g(\mathbf{x}_k)\right) \quad \forall s_0 \in \mathbb{R}.
\end{equation}
Let us estimate the right-hand side of \eqref{k_min_ineq}. By the $L$-smoothness of $f$, for $\mathbf{y}:=\mathbf{x}_k + s_0 g(\mathbf{x}_k)$ and $\mathbf{x}=\mathbf{x}_k$ we have
$$f(\mathbf{y}) \leq f(\mathbf{x}_k) + s_0\langle \nabla f(\mathbf{x}_k), g(\mathbf{x}_k)\rangle + \frac{Ls_0^2}{2}\|g(\mathbf{x}_k)\|^2.$$
Since $\langle \nabla f(\mathbf{x}_k), g(\mathbf{x}_k)\rangle = \|g(\mathbf{x}_k)\|^2 + \langle \nabla f(\mathbf{x}_k) - g(\mathbf{x}_k), g(\mathbf{x}_k)\rangle$ and $|\langle \nabla f(\mathbf{x}_k) - g(\mathbf{x}_k), g(\mathbf{x}_k)\rangle| \leq \delta_1\|g(\mathbf{x}_k)\|$ by \eqref{g_cond_i}, using Young's inequality $|s_0|\delta_1\|g(\mathbf{x}_k)\| \leq \frac{Ls_0^2}{2}\|g(\mathbf{x}_k)\|^2 + \frac{\delta_1^2}{2L}$, we can conclude that
\begin{equation}
    f(\mathbf{x}_{k+1}) \leq f(\mathbf{x}_k) +\left(s_0 + s_0^2 L\right)\|g(\mathbf{x}_k)\|^2 + \frac{1}{2L}\delta_1^2
\end{equation}
for each $s_0\in\mathbb{R}$. Further,
\begin{equation}\label{eqineqs0}
     -\left(s_0 + s_0^2 L\right) \|g(\mathbf{x}_k)\|^2 \leq f(\mathbf{x}_k) - f(\mathbf{x}_{k+1}) + \frac{1}{2L}\delta_1^2.
\end{equation}
The inequalities above allow one to obtain the result about the convergence of the SESOP method (see \cite{kuruzov_stonyakin}, where the proofs of Theorems \ref{SESOP_theorem_full} and \ref{delta_link} below can be found).
\begin{theorem}
\label{SESOP_theorem_full}
Let the objective function $f$ be $L$-smooth and $\gamma$-quasar-convex with respect to $\mathbf{x}^*$. Let $\tau_k$ be the step value obtained with the inexact solution of the auxiliary problem \eqref{eqsubproblem_sesop} on step 2 in Algorithm \ref{alg:sesop} on the $k$-th iteration. Namely, the following conditions for inexactness hold:
\begin{itemize}
    \item[(i)] For the inexact gradient $g:\mathbb{R}^n\rightarrow\mathbb{R}^n$ there is some constant $\delta_1\geq 0$ such that for all points $\mathbf{x}\in\mathbb{R}^n$ condition \eqref{g_cond_i} holds.
    \item[(ii)] The inexact solution $\tau_k$ meets the following condition:
    \begin{equation}
\label{grad_orth_i}
|\left\langle \nabla f(\mathbf{x}_k), \mathbf{d}_{k-1}^2 \right\rangle |\leq k^2 \delta_2
\end{equation}
for some constant $\delta_2 \geq 0$ and each $k\in\mathbb{N}$. Note that $\mathbf{x}_k=\mathbf{x}_{k-1}+D_{k-1}\tau_{k-1}$.
    \item[(iii)] The inexact solution $\tau_k$ meets the following condition for some constant $\delta_3 \geq 0$:
    \begin{equation}
    \label{grad_orth_ii}
        |\left\langle \nabla f(\mathbf{x}_k), \mathbf{x}_k - \mathbf{x}_0\right\rangle| \leq \delta_3.
    \end{equation}
    \item[(iv)] The problem from step 2 in Algorithm \ref{alg:sesop} is solved with accuracy $\delta_4 \geq 0$ on the function on each iteration, i.e.
    $f(\mathbf{x}_k) - \min_{\tau\in\mathbb{R}^n}f(\mathbf{x}_{k-1} + D_{k-1}\tau) \leq \delta_4$.
\end{itemize}
Then the sequence $\{\mathbf{x}_k\}$ generated by Algorithm \ref{alg:sesop} satisfies
\begin{equation}
 f(\mathbf{x}_k)-f^* \leq \frac{8LR^2}{\gamma^2 k^2} + \left(\frac{R}{\gamma}+ 10\right)\delta_1 + 4\sqrt{\delta_2} +  \delta_3 + 5\sqrt{\frac{L\delta_4}{k}}
\end{equation}
for each $k\geq 8$, where $R=\|\mathbf{x}^*-\mathbf{x}_0\|$.
\end{theorem}

Moreover, conditions (ii) and (iii) are implied by condition (iv) with the appropriate constants $\delta_2, \delta_3$.
\begin{theorem}
\label{delta_link}
If condition $(iv)$ from Theorem \ref{SESOP_theorem_full} holds, then conditions \eqref{grad_orth_i} and \eqref{grad_orth_ii} hold with the constants $\delta_2, \delta_3 \geq 0$ satisfying the following estimates:
$$\delta_3 \leq \sqrt{2L \delta_4}\left(\sqrt{\max_{k}(\|D_k\| \|\tau_k\|)} + \sqrt{\max_{k}\|\mathbf{d}^1_{k-1}\|}\right)$$
and
$$\delta_2 \leq \frac{1}{k^2}\sqrt{2L \max_{k}\|\mathbf{d}_k^2\|\delta_4}.$$
\end{theorem}

\begin{remark}
Let us comment on the verifiability of conditions (ii) and (iii) of Theorem \ref{SESOP_theorem_full}. Although they are formulated in terms of the exact gradient and cannot be checked directly, Theorem \ref{delta_link} shows that they are implied by condition (iv), i.e. by the $\delta_4$-accuracy of the solution of the auxiliary problem in the function value, which is controllable: it can be guaranteed by a certifiable method for the subproblem (e.g., by the ellipsoid method with the corresponding number of iterations, see \eqref{N_ellispoids}). The quantities in the right-hand sides of the estimates of Theorem \ref{delta_link} are bounded whenever the trajectory of the method is bounded, and they are available during the run of the algorithm (see the remark after Theorem \ref{theorem_exact_tau}).
\end{remark}

Using the results above, we are ready to estimate the complexity of the SESOP method.

\section{Auxiliary Low-dimensional Subproblems}

In this section, we provide the analysis of the complexity of the auxiliary low-dimensional subproblems of the SESOP method.

\subsection{Proof of Theorem \ref{theorem_exact_tau}}
\label{proof_theorem_exact_tau}
Theorem \ref{delta_link} gives the correspondence between the inaccuracies $\delta_2, \delta_3$ and $\delta_4$:
$$\delta_3 \leq \sqrt{2L \delta_4}\left(\sqrt{\max_{k}(\|D_k\| \|\tau_k\|)} + \sqrt{\max_{k}\|\mathbf{d}^1_{k-1}\|}\right)$$
and
$$\delta_2 \leq \frac{1}{k^2}\sqrt{2L \max_{k}\|\mathbf{d}_k^2\|\delta_4}.$$
The condition
\begin{equation}
\label{max_epsilon}
\max\left\{\frac{8LR^2}{\gamma^2 k^2}, \left(\frac{R}{\gamma}+ 10\right)\delta_1, 4\sqrt{\delta_2},  \delta_3, 5\sqrt{\frac{L\delta_4}{k}}\right\} \leq \frac{\varepsilon}{5},
\end{equation}
is sufficient to approach the quality $\varepsilon$ in the function value by Theorem \ref{SESOP_theorem_full}. So, according to \eqref{max_epsilon}, we have the following conditions on the inexactnesses $\delta_2, \delta_3, \delta_4$:
\begin{equation}
\label{g_orth_i_eps}
\delta_2 \leq \frac{\varepsilon^2}{400}
\end{equation}

\begin{equation}
\label{g_orth_ii_eps}
        \delta_3 \leq \frac{\varepsilon}{5}
\end{equation}

\begin{equation}
\label{delta_4_eps}
\delta_4\leq \frac{\varepsilon^2}{625L}
\end{equation}

So we obtain the main statement about the quality of the subproblem solution and the conditions on the number of iterations and the inexactness of the gradient.

\begin{lemma}
\label{exact_tau_lemma}
Let the inexact gradient meet condition \eqref{g_cond_i} with $$\delta_1\leq \frac{\varepsilon}{\frac{R}{\gamma}+10}.$$ To obtain the quality $\varepsilon$ in the function value, the SESOP method (Algorithm \ref{alg:sesop}) with the inexact gradient should make $$T \geq \sqrt{\frac{40LR^2}{\gamma^2 \varepsilon}}$$ iterations, and at each iteration the subproblem should be solved so that conditions \eqref{g_orth_i_eps}, \eqref{g_orth_ii_eps} and \eqref{delta_4_eps} are met.
\end{lemma}

Moreover, uniting the conditions \eqref{g_orth_i_eps}--\eqref{delta_4_eps} with the results of Theorem \ref{delta_link}, we can obtain a sufficient quality of the subproblem solution for these conditions to hold.

\begin{lemma}
\label{delta_4_cond}
The following quality of the subproblem solution at the $k$-th iteration is sufficient for conditions \eqref{g_orth_i_eps}, \eqref{g_orth_ii_eps} and \eqref{delta_4_eps} to be met:
$$\delta_4 \leq \min\left\{\frac{\varepsilon^2}{50 L\left(\sqrt{\max_{k}(\|D_k\| \|\tau_k\|)} + \sqrt{\max_{k}\|\mathbf{d}^1_{k-1}\|}\right)^2}, \frac{\varepsilon^4}{320000\, L \max_{k}\|\mathbf{d}_k^2\|}, \frac{\varepsilon^2}{625L}\right\}.$$
\end{lemma}

\begin{proof}
By Theorem \ref{delta_link}, condition \eqref{g_orth_ii_eps} ($\delta_3\leq\varepsilon/5$) holds if $\sqrt{2L\delta_4}\left(\sqrt{\max_{k}(\|D_k\| \|\tau_k\|)} + \sqrt{\max_{k}\|\mathbf{d}^1_{k-1}\|}\right)\leq \frac{\varepsilon}{5}$, which is equivalent to the first term of the minimum. Condition \eqref{g_orth_i_eps} ($\delta_2\leq\varepsilon^2/400$) holds if $\frac{1}{k^2}\sqrt{2L \max_{k}\|\mathbf{d}_k^2\|\delta_4}\leq \frac{\varepsilon^2}{400}$ for all $k\geq 1$; since the left-hand side is largest at $k=1$, it suffices that $\sqrt{2L \max_{k}\|\mathbf{d}_k^2\|\delta_4}\leq \frac{\varepsilon^2}{400}$, which (squaring, using $400^2\cdot 2 = 320000$) is equivalent to the second term. The third term is condition \eqref{delta_4_eps} itself.
\end{proof}

In other words, for $\varepsilon\leq 1$, it is sufficient to solve the auxiliary problem at each iteration with the accuracy $$\delta_4 = \frac{\varepsilon^4}{320000\, C_k L},$$ where the constant $C_k= 1 + \left(\sqrt{\max_{k}(\|D_k\| \|\tau_k\|)} + \sqrt{\max_{k}\|\mathbf{d}^1_{k-1}\|}\right)^2+\max_{k}\|\mathbf{d}_k^2\|$ is defined by the sequence $\{x_k\}$ generated by the algorithm (indeed, $C_k$ dominates each of the three denominators in the lemma above).

This quality corresponds to the worst case at which the subproblem procedure will be stopped. One may expect that for most problems, conditions \eqref{g_orth_i_eps}, \eqref{g_orth_ii_eps} and \eqref{delta_4_eps} will be met significantly earlier.

The statements above are true for both settings (with exact and inexact low-dimensional gradients), and they allow estimating the number of exact and inexact gradient calculations in these cases. As in the statement of Theorem \ref{theorem_exact_tau}, we suppose that on each iteration we have a ball $\mathcal{B}_{R_\tau}^k\subset \mathbb{R}^3$ of radius $R_\tau$ such that $\tau_k \in \mathcal{B}_{R_\tau}^k$. So, for the first setting we have the following estimates. To approach the quality $\delta_4$ above, according to \eqref{N_ellispoids} we need 
$$18 \ln \frac{640000\, L B_k C_k}{\varepsilon^4 }$$
iterations of the ellipsoid method, where $B_k=\max_{\tau\in\mathcal{B}_{R_\tau}^k} f_k(\tau) - \min_{\tau\in\mathcal{B}_{R_\tau}^k} f_k(\tau)$. Multiplying by the number of outer iterations from Lemma \ref{exact_tau_lemma}, we obtain the statement of Theorem \ref{theorem_exact_tau}.

\subsection{Proof of Theorem \ref{ellipsoids_sesop_final_theorem}}
\label{proof_ellipsoids_sesop_final_theorem}
\begin{proof}
Let us estimate the inexactness in the internal problems. By the chain rule, we have the following equality:
$$\frac{d}{d\tau}f(x_k+D_k\tau)= D_k^\top \nabla f(x)\Big|_{x=x_k+D_k\tau}.$$
So, let us consider the ellipsoid method with the inexact gradient in the following form:
$$g_k(\tau)= D_k^\top g(x_k+D_k\tau).$$
For such a gradient, we have the following estimate for the inexactness:
\begin{equation}
\label{inexactness_auxillary}
\|g_k(\tau)-\nabla_{\tau} f_k(\tau)\|\leq \|D_k\|_2 \delta_1.
\end{equation}
So, to approach the quality $\varepsilon$, we need stronger conditions on the inexactness of the gradient: according to \cite{Kuruzov_Gladin}, to approach the quality $\delta_4$ in the subproblem by the ellipsoid method with the $\|D_k\|_2\delta_1$-inexact gradient, it is sufficient that the following condition holds:
$$\delta_1\leq \frac{\delta_4}{\|D_k\|_2 }.$$
Under the assumptions of Theorem \ref{theorem_exact_tau}, the accuracy $$\delta_4=\frac{\varepsilon^4}{320000\, C_k L}$$ is sufficient, where $C_k$ is defined as after Lemma \ref{delta_4_cond}. So, for $\delta_1$ we have the following condition:
$$\delta_1\leq \frac{\varepsilon^4}{320000\, C_k L \|D_k\|_2 }$$
for all $k$, or
$$\delta_1\leq \frac{\varepsilon^4}{320000\, A_N L},$$
where $A_N=C_N \max_{k} \|D_k\|_2.$ The estimates for the numbers of gradient calculations follow as in the proof of Theorem \ref{theorem_exact_tau} (see Appendix \ref{proof_theorem_exact_tau}).
\end{proof}

\subsection{Ellipsoid Method}
The ellipsoid method is a low-dimensional optimization method. It constructs a sequence of ellipsoids containing the solution, and each next ellipsoid has a significantly smaller volume than the previous one. Let us provide the pseudocode of the ellipsoid method from \cite{Kuruzov_Gladin} (see Algorithm \ref{alg:ellispoids}).
\begin{algorithm}
	\caption{Ellipsoids Method with $\delta$-subgradient.}
	\label{alg:ellispoids}
	\begin{algorithmic}[1]
		\REQUIRE Number of iterations $N \geqslant 1$, $\delta \geqslant 0$, ball $\mathcal{B}_{\mathcal{R}} 	\supseteq Q_x$, its center $c$ and radius $\mathcal{R}$.
		\STATE $\mathcal{E}_0 := \mathcal{B}_{\mathcal{R}},\quad H_0 := \mathcal{R}^2 I_n,\quad c_0 := c$.
		\FOR{$k=0,\, \dots, \, N-1$}
		    \IF {$c_k \in Q_x$}
		        \STATE $w_k := w \in \partial_{\delta} g(c_k)$, 
		        \IF {$w_k = 0$}
		            \RETURN $c_k$, 
		        \ENDIF
		    \ELSE
		        \STATE $w_k := w$, where $w \neq 0$ is such that $Q_x \subset \{ x \in \mathcal{E}_k: \langle w, x-c_k \rangle \leqslant 0 \}.$
		    \ENDIF
		    \STATE $c_{k+1} := c_k - \frac{1}{n+1}\frac{H_k w_k}{\sqrt{w_k^T H_k w_k}}$, \\
		    $H_{k+1} := \frac{n^2}{n^2-1} \left( H_k - \frac{2}{n+1}\frac{H_k w_k w_k^T H_k}{w_k^T H_k w_k} \right)$, \\
		    $\mathcal{E}_{k+1} := \{x: (x-c_{k+1})^T H_{k+1}^{-1} (x-c_{k+1}) \leqslant 1 \}$,
		\ENDFOR
		\RETURN $x^N = \arg\min\limits_{x \in \{c_0, ..., c_N \} \cap Q_x } g(x)$.
	\end{algorithmic}
\end{algorithm}

\end{document}